\documentclass[twoside,12pt]{article}
\usepackage{fullpage}
\usepackage{epsf}

\usepackage{latexsym}
\usepackage{amsfonts}
	
	\DeclareMathSymbol{\R}{\mathalpha}{AMSb}{"52}
	\DeclareMathSymbol{\C}{\mathalpha}{AMSb}{"43}
 \DeclareMathSymbol{\Q}{\mathalpha}{AMSb}{"51}
 \DeclareMathSymbol{\IP}{\mathalpha}{AMSb}{"50}
\newcommand{\rs}{\mbox{$\widehat{\C}$}}

\def\FFF{{\mathcal F}}
\def\GGG{{\mathcal G}}

\newtheorem{thm}{Theorem}[section]
\newtheorem{defn}{Definition}[section]

\newtheorem{prop}{Proposition}[section]

\newtheorem{lemma}{Lemma}[section]
\newtheorem{cor}{Corollary}[section]

\newcommand{\qed}{\nopagebreak \hfill \rule{2mm}{2.5mm} \break \par}

\newcommand{\bdry}{\partial}			 
\newcommand{\id}{\mbox{id}}			 
\newcommand{\cl}{\overline}			 
\newcommand{\Aut}{\mbox{\rm Aut}}			 
\newcommand{\union}{\cup}			 

\newcommand{\mtwo}[4]				 
{\mbox{$\left(\begin{array}{cc} 		 
#1 & #2 \\
#3 & #4 
\end{array}
\right)$}}

\newcommand{\pf}{\noindent {\bf Proof: }}

\newcommand{\be}{\begin{enumerate}}
\newcommand{\eb}{\end{enumerate}}
\newcommand{\bi}{\begin{itemize}}
\newcommand{\ib}{\end{itemize}}

\newcommand{\gap}{\vspace{5pt}}			

\newcommand{\GGamma}{\mbox{ $\mbox{Gal}(\overline{\Q}/ \Q)\;$}}
\newcommand{\Gal}{\mbox{\rm Gal}}
\newcommand{\DBP}{\mbox{ $DBP\;$}}
\newcommand{\Qbar}{\mbox{$\overline{\Q}\;$}}
\newcommand{\CP}{\mbox{$\IP^1 \C \;$}}
\newcommand{\Inv}{\mbox{\rm Inv}}
\newcommand{\Stab}{\mbox{\rm Stab}}
\newcommand{\Coeff}{\mbox{{\small Coeff}}}
\newcommand{\Moduli}{\mbox{{\small Moduli}}}
\newcommand{\Rat}{\mbox{\rm Rat}}
\newcommand{\PSL}{\mbox{PSL}}
\newcommand{\Fix}{\mbox{\rm Fix}}
\newcommand{\Mon}{\mbox{\rm Mon}}

\def\IMSmarkvadjust{0 pt}
\def\IMSmarkhadjust{0 pt}
\def\IMSmarkhpadding{0 pt}
\def\SBIMSMark#1#2#3{
 \font\SBF=cmss10 at 10 true pt
 \font\SBI=cmssi10 at 10 true pt
 \setbox0=\hbox{\SBF \hbox to \IMSmarkhpadding{\relax}
                Stony Brook IMS Preprint \##1}
 \setbox2=\hbox to \wd0{\hfil \SBI #2}
 \setbox4=\hbox to \wd0{\hfil \SBI #3}
 \setbox6=\hbox to \wd0{\hss
             \vbox{\hsize=\wd0 \parskip=0pt \baselineskip=10 true pt
                   \copy0 \break%
                   \copy2 \break%
                   \copy4 \break}}
 \dimen0=\ht6   \advance\dimen0 by \vsize \advance\dimen0 by 8 true pt
                \advance\dimen0 by -\pagetotal
	        \advance\dimen0 by \IMSmarkvadjust
 \dimen2=\hsize \advance\dimen2 by .25 true in
	        \advance\dimen2 by \IMSmarkhadjust

  \openin2=publishd.tex
  \ifeof2\setbox0=\hbox to 0pt{}
  \else 
     \setbox0=\hbox to 3.1 true in{
                \vbox to \ht6{\hsize=3 true in \parskip=0pt  \noindent  
                {\SBI Published in modified form:}\hfil\break
                \input publishd.tex 
                \vfill}}
  \fi
  \closein2
  \ht0=0pt \dp0=0pt
 \ht6=0pt \dp6=0pt
 \setbox8=\vbox to \dimen0{\vfill \hbox to \dimen2{\copy0 \hss \copy6}}
 \ht8=0pt \dp8=0pt \wd8=0pt
 \copy8
 \message{*** Stony Brook IMS Preprint #1, #2. #3 ***}
}

\begin{document}

\title{Dessins d'enfants and Hubbard trees}

\author{
		Kevin M. Pilgrim \footnotemark[1]\\ 
		Dept. of Mathematics and Statistics\\
		University of Missouri at Rolla\\
		Rolla, MO 65409-0020\\
		USA\\
  pilgrim@umr.edu	
  }
\renewcommand{\today}{\relax}

\maketitle 
\SBIMSMark{1999/5}{April 1999}{}
\thispagestyle{empty}

\begin{abstract}
We show that the absolute Galois group $\GGamma$ acts faithfully on
the set of Hubbard trees.  Hubbard trees are finite planar trees,
equipped with self-maps, which classify postcritically finite
polynomials as holomorphic dynamical systems on the complex plane.  
We establish an explicit relationship between certain Hubbard trees
and the trees known as ``dessins d'enfant'' introduced by
Grothendieck.
\end{abstract}

\tableofcontents

  \footnotetext[1]{
		Research supported in part by NSF Grant DMS-97003724.  
		{\it Key words and phrases:}  Belyi, dessins, holomorphic dynamics, Julia
set.
		}
  
\newpage

\section{Introduction}

Recently there has been an attempt to gain an understanding of the
structure of the {\em absolute Galois group} $\Gamma=\GGamma$ by
exploiting the remarkable fact that there is a faithful action of
$\Gamma$ on a certain infinite set of finite, planar trees, called {\em
dessins}.  These dessins are combinatorial objects which classify planar
covering spaces $X \stackrel{f}{\to} \C - \{0,1\}$ given by
polynomial maps $f$ unramified above $\{0,1\}$.  The action of $\Gamma$ on
dessins is obtained by letting $\Gamma$ act on the coefficients of $f$,
which one may take to be algebraic.  

The main result of this paper (Theorem 3.8) is
that there is also a faithful action of
$\Gamma$ on the infinite set of {\em Hubbard trees}, which are finite
planar trees equipped with self-maps, and which arise in the study of
holomorphic dynamical systems.   These Hubbard trees are combinatorial
objects which classify {\em postcritically finite} polynomials $f: \C \to
\C$ as dynamical systems (a polynomial $f$ is postcritically finite if
the {\em postcritical set} $P_f = \cup_{n>0} f^{\circ n}(C_f)$ is finite,
where $C_f$ is the set of critical points in $\C$).   Again, one may take
the coefficients of such a map to be algebraic, and the action of
$\Gamma$ is obtained by letting $\Gamma$ act on the coefficients of $f$.  
In fact, we prove that $\Gamma$ acts faithfully on a highly restricted
subset $DBP$ (''dynamical Belyi polynomials'') consisting of
postcritically finite polynomials
$f$ whose iterates are all unramified over $\{0,1\}$ and whose Hubbard
tree is uniquely determined by the dessin associated to $f$ as a
covering space, plus a small amount of additional data (see Definition
3.3).

There are several intriguing aspects to this dynamical point of view.  
First, it turns out that the natural class of objects with which to
work consists of actual polynomials as opposed to equivalence classes
of polynomials.   Second, the dynamical theory is richer.  In
particular, we will introduce a special class of dynamical Belyi
polynomials which we call {\it extra-clean} and which is closed under
composition, hence under iteration.  This will allow us to associate a
{\it tower} of invariants to a single given polynomial $f$, namely the
monodromy groups $\Mon(f^{\circ n})$ of its iterates.  Finally, the
dynamical theory here embeds into the non-dynamical one in the
following sense:  there is a $\Gamma$-equivariant injection of the set
of extra-clean {\it dynamical} Belyi polynomials into the set of {\it
non-dynamical} isomorphism classes of Belyi polynomials given by $f
\mapsto f^{\circ 2}$ (Theorem \ref{thm:xdbp}).  From the point of view
of dynamics, this is remarkable:  the dynamics of such an $f$, which
involves an identification of domain and range, is completely
determined by the isomorphism class of $f^{\circ 2}$ as a {\it
covering space}, which does not require such an identification.

\gap

\noindent{\bf Organization of this paper.}  In \S 2 we recall the
Grothendieck correspondence giving the combinatorial classification of
algebraic curves defined over $\Qbar$; throughout, we concentrate on
the case of polynomials and planar tree dessins.  In \S 3, we
introduce dynamical Belyi polynomials, relate them with non-dynamical
ones via the notion of a normalization, and prove a preliminary
variant (Theorem \ref{thm:main1}) of our main result in terms of
normalized dessins.  In \S 4 we discuss the use of towers of monodromy
groups to distinguish Galois orbits, give some examples, and derive
recursive formulae for monodromy generators of iterates of maps
(Theorem \ref{thm:recursive}).  In \S 5 we discuss various algebraic
invariants, e.g. fields of moduli and definition, attached to
dynamical Belyi polynomials. We prove (Theorem \ref{thm:moduli}) that
the field $K_{\Coeff}(f)$ generated by the coefficients of a
dynamical Belyi polynomial $f$ coincides with the {\it field of
moduli} of the conjugacy class of $f$ introduced by Silverman
\cite{Si}.  \S 6 is essentially independent, consisting of a
translation of our preliminary main theorem into the language of
Hubbard trees.

\gap

\noindent{\bf Acknowledgements.} I thank A. Epstein, X.  Buff, and
J. Hubbard for useful conversations, C. Henriksen for assistance
drawing Julia sets.  I
am grateful to W. Fuchs for planting the seed from which this project
grew by asking me explicitly some years ago if there was a connection
between dessins and Hubbard trees.
\section{Dessins d'Enfant}

In 1979, Belyi proved a remarkable theorem: an algebraic
curve $X$ defined over $\C$ is defined over $\Qbar$ only if there is a
holomorphic function $f: X \to \CP$, called a {\it Belyi morphism},
all of whose critical values lie in $\{0, 1, \infty\}$, i.e. $X$ is a
branched covering of $\CP$ ramified only over $\{0,1,\infty\}$
\cite{Bel}.  This in turn led Grothendieck to the observation that
there is a faithful action of the absolute Galois group $\Gamma =
\GGamma$ on a set of simple, concrete, combinatorial objects, called
{\it dessins}, which in fact one may take to be certain finite planar
trees.  The structure of the orbits of $\Gamma$ under this action
remains quite mysterious, and the development of effective
combinatorial invariants for distinguishing them has been the subject
of recent work; see e.g. \cite{JS}.
\gap

\noindent{\bf Combinatorial classification of algebraic curves.}  We
begin by outlining the combinatorial classification of algebraic
curves $X$ defined over $\Qbar$.  For a good introduction to the
subject, see e.g. the article by Schneps in \cite{Sch}.  The
statements are cleanest provided one first introduces a minor,
commonly adopted technical notion.  A Belyi morphism $f: X \to \CP$ is
called {\it clean} if the ramification at each point lying over 1 is
exactly equal to two.  Let $q(z) = 4z(1-z)$.  If $f$ is a Belyi
morphism, then $q\circ f = 4f(1-f)$ is a clean Belyi morphism, and so
$X$ is defined over $\Qbar$ if and only if there is a clean Belyi
morphism from $X$ to $\CP$.  If $f$ is clean we call the pair $(X, f)$
a {\it clean Belyi pair}.  Two such pairs $(X_1, f_1), (X_2, f_2)$ are
called {\it isomorphic} if there is an isomorphism $\phi: X_1 \to X_2$
with
$f_1 = f_2 \circ \phi$.  We will be mainly interested in the case when
$X=\CP$ and $f$ is a polynomial.

On the combinatorial side, a {\it Grothendieck dessin} (``scribble''?)
is an abstract simplicial 2-complex with $0$-cells, $1$-cells, and
$2$-cells denoted respectively $X_0, X_1, X_2$ such that the
underlying space is homeomorphic to a closed, connected, oriented
surface, and such that the vertices $X_0$ are given a bipartite
structure, i.e. can be colored black and white such that each edge of
$X_1$ has exactly one black and one white vertex.  Two such dessins
are called {\it isomorphic} if there is an orientation- and
color-preserving isomorphism of complexes.  A dessin is called {\it
clean} if each white vertex is the endpoint of exactly two edges.  The
{\it genus} of a dessin is the genus of the underlying surface.  We
will be mainly interested in the case when the genus is zero and the
union of edges and vertices forms a tree, i.e. there is a single
two-cell in $X_2$.  In this case we shall specify a dessin by
specifying a planar tree with a bicoloring of vertices.  

A clean Belyi pair $(X, f)$ determines a clean dessin ${\bf D}_f$
whose white vertices are preimages of 1, whose black vertices are
preimages of 0, whose edges are preimages of the segment $[0,1]$, and
whose $2$-cells are preimages of $\CP - [0,1]$.  The classification
may now be formulated as follows (\cite{Sch}, Thm. I.5):  

\begin{thm}{\bf (Grothendieck correspondence) }
\label{thm:correspondence}
The map $(X,f) \mapsto {\bf D}_f$ descends to a bijection between
isomorphism classes of clean Belyi pairs and clean dessins.
\end{thm}

\noindent{\bf Convention.}  We are mainly concerned with
the case when $X = \CP =\rs =\C \union \{\infty\}$ and $f$ is a polynomial.  A
clean Belyi pair $(X,f)$ is determined by a clean Belyi polynomial $f
\in \Qbar [z]$.  Two clean Belyi polynomials are then isomorphic if
and only if there is an affine map $A \in \Aut (\C)$ with $f = g\circ
A$.  Note that if $f = g \circ A$, with $f, g \in \Qbar [z]$ of degree
at least one, then necessarily $A \in \Qbar [z]$ since $A$ must send
the set $f^{-1}(\{0,1\})$ onto the set $g^{-1} (\{0,1\})$ and both
sets consist of a collection of at least two algebraic numbers.

Throughout the remainder of this work, {\it we will deal exclusively with
clean Belyi polynomials and clean dessins of genus 0, i.e. planar tree
dessins.  We therefore now adopt the convention that the term
``dessin'' means clean planar tree dessin, and that ``Belyi polynomial''
means clean Belyi polynomial, unless otherwise specified.  }
\gap

\noindent{\bf Notation.}  
\begin{itemize}
\item $\Aut (\C)$, the group of affine maps $az+b, a \neq 0$;
\item $A,B$, elements of $\Aut (\C)$;
\item $f,g, f_1, f_2$, clean Belyi polynomials;
\item $BP$, the set of clean Belyi polynomials;
\item $[BP]$, the set of isomorphism classes of clean Belyi
polynomials;
\item $[f]$, the isomorphism class of $f$ as a Belyi polynomial;
\item ${\bf D} _f$, the dessins of $f$;
\item $[{\bf D} _f]$, the isomorphism class of ${\bf D} _f$;
\item $\Gamma$, the absolute Galois group $\GGamma$.  
\end{itemize}

\noindent{\bf Action of \GGamma on dessins.}  
The group $\Gamma$ acts on clean Belyi
polynomials by twisting coefficients, i.e. if $\sigma \in
\Gamma$ and $f(z) = a_d z^d + ... + a_0$ then 
\[ f^{\sigma} = \sigma(a_d)z^d + ... +
\sigma(z_0).\] 
This action descends to an action on the set $[BP]$ of
isomorphism classes of Belyi polynomials, since if $f_1 = f_2 \circ
A$, then $f_1 ^{\sigma} = (f_2 \circ A)^{\sigma} = f_2 ^{\sigma} \circ
A^{\sigma}$.  By the Grothendieck correspondence, we get an action of
$\Gamma$ on isomorphism classes of dessins.  Lenstra and Schneps
(\cite{Sch}, Thm. II.4) have shown

\begin{thm}
\label{thm:faithful}
The action of $\Gamma$ on the set $[BP]$ of isomorphism
classes of Belyi polynomials, hence on the set of 
dessins, is faithful.  
\end{thm}

In fact, their argument is constructive: given any $\sigma \in \Gamma$
and $\alpha, \beta \in \Qbar$ with $\sigma(\alpha) = \beta \neq
\alpha$, they produce, by using the arguments in the proof of Belyi's
theorem and an elementary, technical, algebraic lemma (Lemma
\ref{lemma:lenstra-schneps} below), a pair $f_{\alpha}, f_{\beta}$ of
nonisomorphic Belyi polynomials with $f_{\alpha} ^{\sigma} =
f_{\beta}$.

Since the notion of isomorphism between Belyi polynomials involves a
coordinate change in the domain, but not in the range, they cannot be
considered as dynamical objects.  In the next section, we replace the
notion of isomorphism with that of affine conjugacy, and show that the
action of $\Gamma$ on a suitable set of affine conjugacy classes is
faithful.

\section{Dynamical Belyi polynomials}

Let $f: \C \to \C$ be a polynomial of degree $d \geq 2$.  Complex
dynamics is concerned with the behavior of points under {\it
iteration} of such a function, i.e. with the behavior of {\it orbits}
\[ \{z, f(z), f^{\circ 2}(z), f^{\circ 3}(z), ... \}\] 
where $f^{\circ n}$ denotes the $n$-fold composition of $f$ with
itself.  It turns out that understanding the orbits of the critical
points (i.e. those $c \in \C$ for which $f'(c)=0$) is crucial to
understanding the global dynamics of $f$.  Let $C_f$ denote the set
of critical points of $f$ and $V_f = f(C_f)$ the set of critical
values.  We define the {\it postcritical set} of $f$ by 
\[ P_f = \cl{ \bigcup_{n>0}f^{\circ n}(C_f)}.\]
Then $V_f \subset P_f$, $f(P_f) \subset P_f$, and $P_{f^{\circ n}} =
P_f$ for all $n > 0$.

\begin{defn}A {\rm \bf  dynamical Belyi polynomial} is a 
Belyi polynomial $f$ of degree $d \geq 3$ for which $P_f \subset
\{0,1\}$.  We denote by $\DBP$ the set of all dynamical Belyi
polynomials.
\end{defn}

Recall that, by convention, $f$ is assumed clean.

\begin{prop}Let $f$ be a dynamical Belyi polynomial.  Then 
$V_f = P_f = \{0,1\}$.  Moreover, $f^{-1}(\{0,1\}) \supset \{0,1\}$
and hence $0$ and $1$ are vertices of the dessin ${\bf D}_f$.
\end{prop}

\pf $P_f \subset \{0,1\}$ by definition.  Since $f$ is clean, the
branching above 1 is exactly equal to two, hence $1 \in V_f \subset
P_f$.  There are $d-1$ critical points, counted with multiplicity, and
$d/2$ of them are the preimages of 1, since $f$ is clean.  If $d>2$,
then $(d-1) - d/2 > 0$, so there are critical points which do not map
to 1.  Since $f$ is a Belyi polynomial, the other critical points map
to $0$, so $0 \in V_f \subset P_f$.  The last statement follows since
$f(\{0,1\}) \subset \{0,1\}$.  
\qed

\begin{prop}
\label{prop:Bisid}
Let $f,g \in DBP$, and suppose $g = B \circ f \circ A$,
where $A, B \in \Aut(\C)$.  Then $B=\id$.
\end{prop}

\pf We have 
\[ \{0,1\}= V_g = V_{B\circ f \circ A} = B(V_{f \circ A}) = B(V_f) =
B(\{0,1\})\] 
where the first equality follows from the preceding proposition.  
So $B=\id$ or $B(z)=1-z$.   To rule out the latter case, choose $z
\in f^{-1}(0)$ mapping to $0$ with local degree one.  Then 
\[g(A^{-1}(z)) = (B\circ f \circ A)(A^{-1}(z)) = B(0) = 1\]
mapping by local degree one, which violates the cleanness of $g$.
\qed

Two polynomials $f, g: \C \to \C$ for which there is an affine map $A:
\C \to \C$ satisfying $g = A^{-1}fA$ are called {\it conjugate}.  As
dynamical systems, they are the same, just viewed in different
coordinates.  A conjugacy from $f$ to itself is called an {\it
automorphism} of $f$.  As a corollary to the previous proposition,
upon setting $B=A^{-1}$ we obtain

\begin{thm}
\label{thm:nonconj}
No two distinct elements $f, g$ of $\DBP$ are conjugate,
and no element of $\DBP$ has a nontrivial automorphism.
\end{thm}

\noindent{\bf Remark and convention:} This theorem fails without the
restriction $d \geq 3$, as the maps $f_a (z) = a(z-1)^2 + 1, \; a \in
\C-\{0\} $ are clean Belyi polynomials with $P_{f_a} = \{ 1 \} $
conjugate to $z^2$.  Thus {\it we now assume throughout that the
degrees of all Belyi polynomials are at least three}.  
\gap

\begin{thm}
\label{thm:zofix}
Let $f \in DBP$, and let 
\[ Z_f = \{z | f(z)=0\},\;\; O_f = \{z | f(z)=1\}, \;\; \Fix(f) = \{z\;\;
| \;\;f(z) = z\}.\]
Then $f$ is uniquely determined by any one of the three sets $Z_f,
O_f, \Fix(f)$, counted with multiplicity. 
\end{thm}

\pf If $f,g \in DBP$ and either $Z_f = Z_g$ or $O_f = O_g$, then $f =
B \circ g$ where $B \in \Aut(\C)$, and so by Theorem \ref{thm:nonconj}
$f=g$.  If $\Fix(f) = \Fix(g)$ then since $f(\{0,1\}), g(\{0,1\})
\subset \{0,1\}$ we must have $f|\{0,1\} = g| \{0,1\}$.  Since
$\Fix(f)=\Fix(g)$,
\begin{equation}
 f(z) - z = \lambda(g(z)-z), \;\;\mbox{some} \;\; \lambda \in \C. 
\end{equation}
If either $f(0)=g(0)=1$ or $f(1)=g(1)=0$ then substituting into (1)
implies $\lambda=1$ and $f=g$.  Otherwise, $f(0)=g(0)=0$ and
$f(1)=g(1)=1$.  Differentiating Equation (1) we obtain
\begin{equation}
f'(z)-1 = \lambda(g'(z)-1).
\end{equation}
The cleanness criterion implies that $f'(1)=g'(1)=0$, and substituting
this into Equation (2) implies $\lambda=1$ and $f=g$.
\qed

\noindent{\bf DBPs and normalized Belyi polynomials.}  We next relate
dynamical Belyi polynomials and non-dynamical ones via the notion of a
normalization.

\begin{defn}
A {\bf\rm normalized} Belyi polynomial is a pair 
$(f, (z, w))$ where $f$ is a Belyi polynomial and $(z,w)$ is an
ordered pair of distinct (algebraic) numbers with $\{z, w\} \subset
f^{-1}(\{0,1\})$.  Two normalized Belyi polynomials $(f_1, (z_1
, w_1 ))$ and $(f_2, (z_2, w_2))$ are called {\bf \rm isomorphic} if
there is an $A \in \Aut(\C)$ for which $f_1 = f_2 \circ A$,
$z_2 = A(z_1)$, and $w_2 = A(w_1)$.  We denote the set of isomorphism
classes of normalized Belyi polynomials by $[BP^*]$.  
\end{defn}

Note that we do {\it not} require that $z \in f^{-1}(0)$ and $w \in
f^{-1}(1)$.

\begin{defn}
A {\bf \rm normalized dessin} ${\bf D}^*$ is a dessin ${\bf D}$
together with an ordered pair $(z,w)$ of vertices of ${\bf D}$.  Two
normalized dessins ${\bf D} ^* _1, {\bf D} ^* _2 $ are called {\bf \rm
isomorphic} if there is an isomorphism ${\bf D}_1 \to {\bf D}_2$ of
dessins carrying one ordered pair of vertices to the other.
\end{defn}

An immediate consequence of the definitions and the Grothendieck
correspondence is that the natural map sending a normalized Belyi
polynomial $(f, (z,w))$ to the normalized dessins $({\bf D} _f, (z,w))$
descends to a bijection between isomorphism classes of geometric
objects (normalized Belyi polynomials) and combinatorial ones
(normalized abstract dessins).  

Recall that if $f \in DBP$ then $0$ and $1$ are vertices of ${\bf D}_f$.  
There are natural maps 
\[ 
\begin{array}{lcl}
DBP \to [BP^*] & \mbox{given by } & f \mapsto [(f,(0,1))] \\
DBP \to [BP]   & \mbox{given by } & f \mapsto [f] \\
\mbox{$[BP^* ]$} \to [BP] & \mbox{induced by } & (f, (z,w)) \mapsto [f] \\
\end{array}
\]
where $[f]$ is the isomorphism class of $f$ as a Belyi polynomial.
The following diagram then commutes:

\begin{center}
\begin{picture}(150,150)
\put(75,25){\makebox(0,0){$[BP]$}}
\put(25,125){\makebox(0,0){$DBP$}}
\put(125,125){\makebox(0,0){$[BP^*]$}}
\put(50,125){\vector(1,0){50}}
\put(25,115){\vector(1,-2){37}}
\put(125,115){\vector(-1,-2){37}}
\end{picture}
\end{center}

\begin{thm} The map $DBP \to [BP^*]$ given by 
\[  f \mapsto [(f, (0,1))]  \]
is a bijection.
\end{thm}

\pf The map is clearly injective, since if $(f_1, (0,1))$ is isomorphic
to $(f_2, (0,1))$, then the affine map $A$ giving the isomorphism must
send zero to zero and one to one.  Hence $A$ is the identity and
$f=g$.  The map is surjective as well.  First, any isomorphism class
of normalized clean Belyi polynomial contains a 
representative where $z = 0$ and $w = 1$, which can be constructed as
follows.  Choose any representative $(f, (z, w))$, and let $A$ be the
unique affine map which sends $0$ to $z$ and $1$ to $w$.  Then $(f
\circ A, (0, 1))$ is equivalent to $(f, (z,w))$.
We now claim that $f \circ A \in DBP$.  First, $V_{f \circ A} = V_f =
\{0,1\}$ since $f$ is assumed clean and precomposing $f$ by an affine
map does not change the set of critical values.  On the other hand, by
construction, $f \circ A (\{0,1\}) \subset \{0,1\}$.  Hence $P_{f
\circ A} \subset \{0,1\}$ and so $f \circ A \in DBP$.  By the
definition of isomorphism in $BP^*$, $[f \circ A, (0,1)] = [f, (z,w)]$
and so the map is surjective.
\qed

Thus, a normalized Belyi polynomial $(f,(z,w))$ determines a
holomorphic dynamical system $g = f \circ A \in DBP$ by identifying
range and domain via an affine map $A$ sending $0$ to $z$ and $1$ to
$w$.  The above theorem implies that $DBP$ is in bijective
correspondence with $[BP^*]$, which in turn is in bijective
correspondence with the set of normalized clean dessins.

\gap

\noindent{\bf The fibers of the map $DBP \to [BP]$.}  The fiber of the
forgetful map $DBP \to [BP]$ over a given element $[f] \in [BP]$ is a
disjoint union of four nonempty subsets, which we describe in terms of
the identification of $DBP$ with $[BP^*]$ given above.  In our
normalization of $f$, we may freely and independently choose $z$ or
$w$ to be a black vertex (a preimage of $0$) or a white vertex (a
preimage of $1$) of ${\bf D} _f$, giving us four possibilities, all of
which can occur.  In terms of dynamics, suppose $g$ is the element of
$DBP$ corresponding to $[(f, (z,w))] \in [BP^*$] in Theorem 3.3 .  Then
\[
\begin{array}{cccc}
 z & \mbox{is black} &\iff &g(0)=0\\
 z & \mbox{is white} &\iff &g(0)=1\\
 w & \mbox{is black} &\iff &g(1)=0\\
 w & \mbox{is white} &\iff &g(1)=1
\end{array}
\] 
Within each of these classes, one can further classify points $g$ in the
fiber over $[f]$ by recording the local degrees of $g$ near $0$ and
$1$.  

\gap

\noindent{\bf Extra-clean dynamical Belyi polynomials.} Dynamics is
concerned with iteration, and although an iterate of a dynamical Belyi
polynomial is again a Belyi polynomial, the property of cleanness may
be lost.  For example, if we choose $w$ to be a white vertex, then
$g(1) = 1$ by local degree two, and so $g^{\circ 2} (1) = 1$ but now
mapping by local degree four, violating cleanness.  To remedy this, we
formulate 

\begin{defn} An element $g \in DBP$ is called {\rm extra-clean} if
$g(1)=g(0) = 0$, and the local degrees of $g$ near $0$ and $1$ are both
equal to one.  The set of all extra-clean dynamical Belyi polynomials
we denote by $XDBP$.  
\end{defn}

In terms of normalizations, suppose $g \in DBP$ corresponds to the
class of normalized dessins $[({\bf D}, (z,w))]$.  Then $g \in XDBP$ if
and only if $z,w$ are both ends of the dessins ${\bf D}$, which are
necessarily black since, by cleanness, white vertices are always
incident to two edges.  We denote by $[XBP^*]$ the subset of $[BP^*]$
corresponding to $XDBP$, and refer to the associated normalized
dessins as {\it extra-clean normalized dessins}.  We obtain the 
the following commutative diagram:

\begin{center}
\begin{picture}(150,150)
\put(75,25){\makebox(0,0){$[BP]$}}
\put(25,125){\makebox(0,0){$XDBP$}}
\put(125,125){\makebox(0,0){$[XBP^*]$}}
\put(50,125){\vector(1,0){50}}
\put(25,115){\vector(1,-2){37}}
\put(125,115){\vector(-1,-2){37}}
\end{picture}
\end{center}
where the map $XDBP \to [BP]$ is surjective.  

\gap

\noindent{\bf Remark:}  We emphasize here that $DBP$ is 
a set of maps, not a set of maps modulo an equivalence relation.
Composition does not descend from $BP$ to a well-defined operation
on isomorphism classes of non-dynamical Belyi polynomials; see the
example in Section 4 for two polynomials $f,g$ with ${\bf D}_f$
isomorphic to ${\bf D}_g$ but with ${\bf D}_{g^{\circ 2}}$ and ${\bf
D}_{f^{\circ 2}}$ non-isomorphic.
\gap

Indeed, if $f, g \in XDBP$, this is always the case:

\begin{thm}
\label{thm:xdbp}
If $f, g \in XDBP$, then $f\circ g \in XDBP$.  Moreover,  
\[ {\bf D}_{f^{\circ 2}} \simeq {\bf D}_{g^{\circ 2}} \iff f=g.\]
Thus the map 
\[ XDBP \to [BP] \]
given by 
\[ f \mapsto [f^{\circ 2}] \]
is injective and $\Gamma$-equivariant.
\end{thm}

This is perhaps remarkable, since it implies that the dynamical system
generated by $f$ is completely determined by the topology of $f^{\circ
2}$ as a covering space.   This property fails even for the highly
restricted set of postcritically finite quadratic polynomials $p$. 
Apart from $p(z)=z^2$, $p^{\circ 2}$ will have two finite critical
values, and as covering spaces of the twice-punctured plane the
second iterates of any two such $p$ are isomorphic.  

The proof relies on the following lemma of Lenstra and Schneps used in
their proof of Theorem \ref{thm:faithful}, and a fact from holomorphic
dynamics:

\begin{lemma}{\rm(\cite{Sch}, Lemma II.3)}
\label{lemma:lenstra-schneps}
Suppose $G, H, \tilde{G}, \tilde{H}$ are polynomials with\\
$G \circ H = \tilde{G} \circ \tilde{H}$ and $\deg(H) =
\deg(\tilde{H})$.  Then there exist constants $c, d$ for which
$\tilde{H} = cH+d$, i.e. $\tilde{H} = B \circ H$ for some $B \in
\Aut(\C)$.  
\end{lemma}

\noindent{\bf Proof of Theorem.}  For a polynomial $f$, set $\Fix(f) =
\{p \;\;|\;\; f(p)-p=0\}$.  From dynamics, one knows that if there is
a $p \in \Fix(f)$ which is a multiple root of $f(z)-z$, (i.e. the {\it
multiplier} of $f$ at $p$ is equal to one) then there is a critical
point of $f$ whose forward orbit is infinite and accumlates at $p$
(\cite{Bea}, Theorem 9.3.2).  For dynamical Belyi polynomials this
cannot occur, since $P_f = \{0,1\}$.  Hence, for any dynamical Belyi
polynomial, the fixed points are all simple.  Thus
\[ f(z)-z\;\; = \;\;c \cdot \prod_{p \in \Fix(f)}(z-p)\]
for a nonzero constant $c$.  By Theorem \ref{thm:zofix} it suffices to show
$\Fix(f) = \Fix(g)$.

The hypothesis and the Grothendieck correspondence imply that there is
an $A \in \Aut(\C)$ with $f^{\circ 2} = g^{\circ 2} \circ A$.  We write
this as $f\circ f = g \circ (g \circ A)$.  The polynomials $f$ and $g$
have the same degree, since their second iterates have the same
degree.  Applying Lemma \ref{lemma:lenstra-schneps} with $G=g, H=g\circ A,
\tilde{G}=\tilde{H}=f$ we obtain an affine map $B$ for which $f = B
\circ (g \circ A)$.  By Proposition \ref{prop:Bisid}, $B=\id$ and so $f=g\circ A$.

Now let $p \in \Fix(f)$.  Then by the previous paragraph 
\[f(p) 		 =  g\circ A (p) 		 =  p \]
and by hypothesis
\[ f^{\circ 2}(p)   =  g^{\circ 2} \circ A (p) 	 =  p.\]
Applying $g$ to both sides of the last equality in the first equation, 
and comparing with the second we get 
\[ g(p) = g^{\circ 2} \circ A (p) = p.\]
Hence $\Fix (f) \subset \Fix(g)$.  Equality follows, either by
appealing to the symmetry of the roles of $f$ and $g$, or the fact
that $\Fix(f), \Fix(g)$ have the same size.
\qed
\gap

\noindent{\bf Galois action on $DBP$.}
The group $\Gamma$ acts on polynomials $f \in \Qbar[z]$ by acting on
its coefficients.  Note that this is a left action, i.e. $f^{\sigma
\tau} = (f^{\tau})^{\sigma}$ since e.g. if $z$ is algebraic, 
\[ f^{\sigma} (z) = \sigma \circ f \circ \sigma^{-1} (z).\]
Hence the action of $\Gamma$ on polynomials in $\Qbar[z]$ by twisting
coefficients is natural with respect to the dynamics in the following
sense: if $f$ is defined over $\Qbar$ and $z \in \Qbar$, then
$\sigma(f(z)) = f^{\sigma}(\sigma(z))$.  

Using this, and the fact that the property of being a critical point
is algebraic, it is easy to show that the group action of $\Gamma$ on
$\Qbar[z]$ preserves the set $DBP$.  Similarly, the action of $\Gamma$
must preserve local degrees, i.e. if $f$ maps $x$ to $y$ by local
degree $k$, then $f^{\sigma}$ maps $\sigma(x)$ to $\sigma(y)$ by local
degree $k$.  Hence $\Gamma$ acts on the set $XDBP$ as well.

Recall that $[XBP^*]$ corresponds to the set $XDBP$ under the
bijection given in Figure 2 and Theorem 3.3.  The group $\Gamma$ also
acts on $[XBP^*]$ in the obvious way: \[ \sigma . [(f, (z,w))] =
[f^{\sigma}, (\sigma(z), \sigma(w))].\] Together with the usual action
of $\Gamma$ on $[BP]$ we find as a consequence of the definitions

\begin{thm}
The diagrams in Figures 1 and 2 are equivariant with respect to the
action of $\Gamma$.
\end{thm}

By the theorem of Lenstra and Schneps, the action of $\Gamma$ on
$[BP]$ is faithful.  By the preceding theorem, Galois orbits in $XDBP$
and in $DBP$ lie over orbits in $[BP]$.  But the forgetful maps $XDBP \to [BP]$
and $DBP$ to $[BP]$ are surjective, so we obtain 

\begin{thm}
The actions of $\GGamma$ on $XDBP$ and on $DBP$ by twisting
coefficients are faithful. 
\end{thm}

Since the correspondence between $XDBP$ and $[XBP^*]$ is a bijection,
and the latter set is in bijective correspondence with isomorphism
classes of normalized dessins where the chosen points are both ends of
the dessins, we obtain a combinatorial version of the preceding theorem:

\begin{thm}
\label{thm:main1}
$\GGamma$ acts faithfully on both the set of isomorphism classes of
extra-clean normalized dessins and on the set of normalized (clean)
dessins.   
\end{thm}

In Section 6, we will show that the sets $DBP$ and $XDBP$ are {\it naturally}
isomorphic respectively to respectively the sets $BHT$ of (isomorphism classes
of) clean Belyi-type Hubbard trees and $XBHT$ of extra-clean Belyi type Hubbard
trees. Our main theorem then follows:

\begin{thm}
\label{thm:main2}
$\GGamma$ acts faithfully on both $BHT$ and $XBHT$.  
\end{thm}

\section{Distinguishing Galois orbits}

It is known (see e.g. \cite{JS}) that if $f, g$ are two Belyi
polynomials (not necessarily clean) which are Galois conjugate, and if
${\bf D}_f, {\bf D}_g$ are their corresponding dessins, then

\begin{enumerate}

\item ${\bf D}_f, {\bf D}_g$ have the same number of edges (i.e. $f$
and $g$ have  the same degree);

\item ${\bf D}_f, {\bf D}_g$ have the same set of valencies (i.e. the unordered
sets of local degrees of $f$ and $g$ at points lying over $0$ and $1$
are the same);

\item the {\it monodromy groups} $\Mon(f), \Mon(g)$ are
permutation-isomorphic.  
\end{enumerate}

The monodromy group\footnote{If $F$ is a Belyi polynomial which is
not necessarily clean, the monodromy group of the clean polynomial $q
\circ F$ is usually called the {\it cartographic group} of $F$.}
$\Mon(f)$ of a Belyi polynomial is the monodromy group of the covering
$f: \C - f^{-1}(\{0,1\}) \to \C - \{0,1\}$ and can be described in
many equivalent ways; see e.g. \cite{JS}.  Here, we think of it as a
finite permutation group acting transitively (but not freely) on the
set of edges $E_f$ of ${\bf D}_f$ (equivalently, on the fiber
$f^{-n}(b)$ of an arbitrary basepoint $b \neq 0,1$, usually $b=1/2$),
and define it as the group generated by the two permutations $\sigma_0
(f)$ which rotates edges counterclockwise about black vertices, and
$\sigma_1 (f)$ which rotates edges counterclockwise about white
vertices; see Figure 7.  The polynomial $f$ is clean exactly when
$\sigma_1 (f)$ is a fixed-point free involution.  The condition that
$\Mon(f)$ and $\Mon(g)$ are permutation-isomorphic means that there is
a bijection $\tau: E_f \to E_g$ such that $\Mon(g) = \tau \Mon(f)
\tau^{-1}$.  The dessins ${\bf D}_f$ and ${\bf D}_g$ are isomorphic
exactly when one can choose a single bijection $\tau$ so that
$\sigma_0 (g) = \tau \sigma_0 (f) \tau^{-1}$ and $\sigma_1 (g) = \tau
\sigma_1 (f) \tau^{-1}$ simultaneously; indeed, this is a special case
of the well-known {\it Hurwitz classification} of coverings.

Thus the monodromy group is an invariant of the Galois orbit.
However, even for clean polynomials, it is known that this invariant
can fail to distinguish Galois orbits (see \cite{JS}, Example 6). 
Below, we take up similar considerations for dynamical Belyi
polynomials.

Recall that $\Gamma$ acts faithfully on the set $XDBP$ of extra-clean
dynamical Belyi polynomials, and that $XDBP$ is closed under
composition and iteration.  For any two polynomials $f, g \in \Qbar
[z]$, we have $(f\circ g)^{\sigma} = f^{\sigma} \circ g^{\sigma}$.  So
in particular if $f,g \in XDBP$ and $g=f^{\sigma}$, then for all
$n>0$, $g^{\circ n} = (f^{\circ n})^{\sigma}$.  Hence by (3) above we
have

\begin{thm}
\label{thm:monodromy-tower}
If $f,g \in XDBP$ and $g=f^{\sigma}$, then $\Mon(f^{\circ n})$ and
$\Mon(g^{\circ n})$ are permutation-isomorphic for all $n>0$.
\end{thm}

\noindent{\bf Question: }{\it To what extent is the converse to Theorem
\ref{thm:monodromy-tower} true?} 
\gap

Recall that by Theorem \ref{thm:xdbp}, the dessins of $f^{\circ 2}$
and $g^{\circ 2}$ are distinct if $f \neq g$.  So a tower of
permutation-isomorphic monodromy groups cannot arise from the trivial
situation where the dessins of $f^{\circ n}$ and $g^{\circ n}$ are 
isomorphic for all $n$.  Below, we give a simple example where we use
this criterion to distinguish orbits lying over a single element of
$[BP]$.  \gap

\noindent{\bf An example.}  Figures 1 and 2 show two normalized
dessins of degree eight; the underlying dessin, without normalization,
is that arising from the Belyi polynomial $h(z)=4z^4(1-z^4)=q\circ (z
\mapsto z^4)$ where $q(z) = 4z(1-z)$.  Denoting by $f$ and $g$ the
corresponding dynamical Belyi polynomials, we see that $f = h \circ
A_f, g=h\circ A_g$ where $A_f, A_g$ are affine maps sending the pair
$(0,1)$ to the indicated vertices, which are $(-1,1)$ for $f$ and
$(-1, i)$ for $g$.  Thus $A_f (z) = 2z-1$ and $A_g (z) = (1+i)z-1$.
Since $q$ and $h$ are both defined over $\Q$, it follows that since
$A_f$ and $A_g$ are not Galois conjugate, the maps $f$ and $g$ are not
Galois conjugate either.

\begin{figure}
\begin{center}
\begin{minipage}{.4\hsize}
\label{fig:normdessf.ps}
\centerline{%
\epsfysize=.95\hsize
\epsfbox{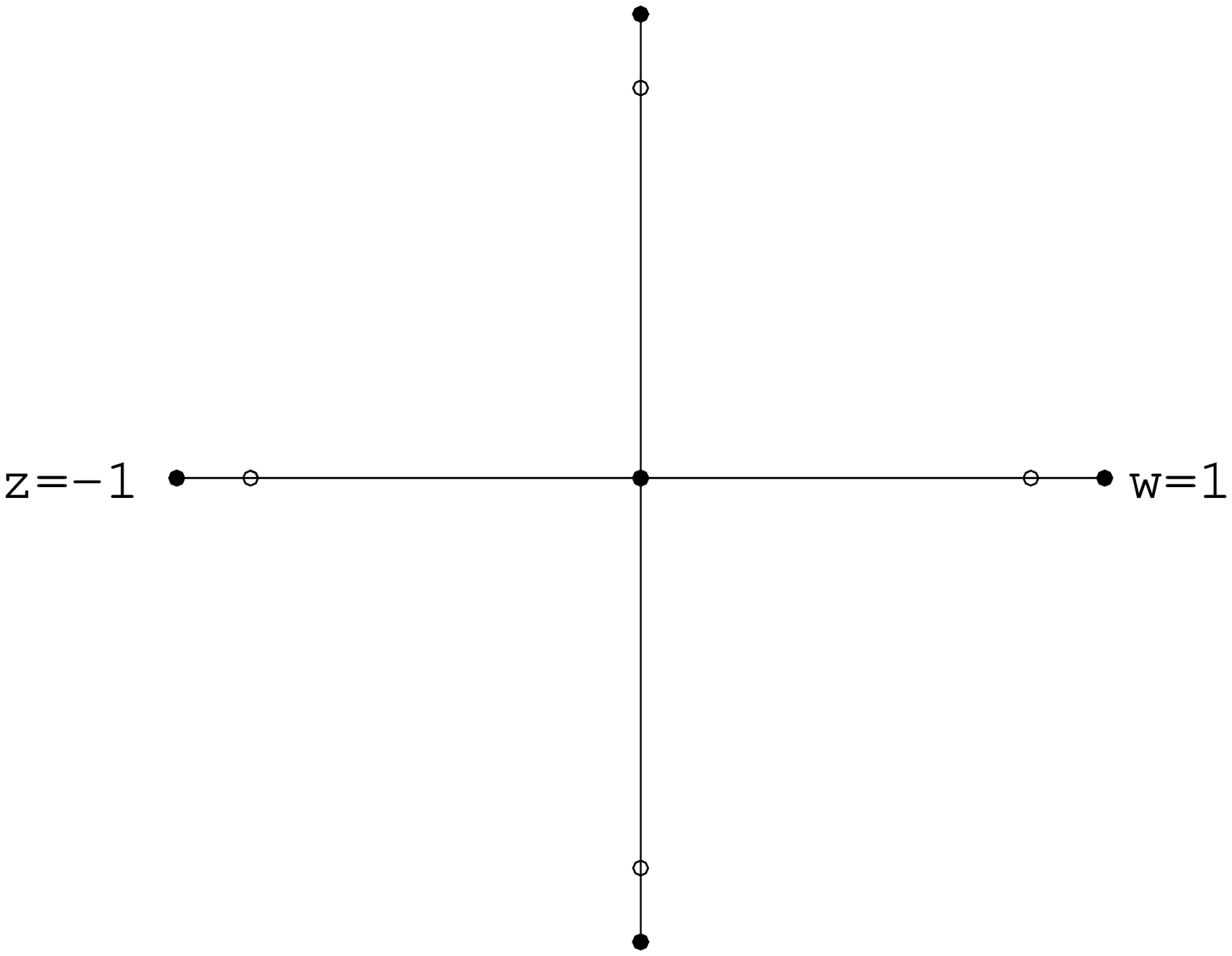}
}
\begin{center}
{\sf Figure 1\\ The normalized dessins of $f$.}
\end{center}
\end{minipage}
%
\qquad\qquad
\begin{minipage}{.4\hsize}
\label{fig:normdessg.ps}
\centerline{%
\epsfysize=.95\hsize
\epsfbox{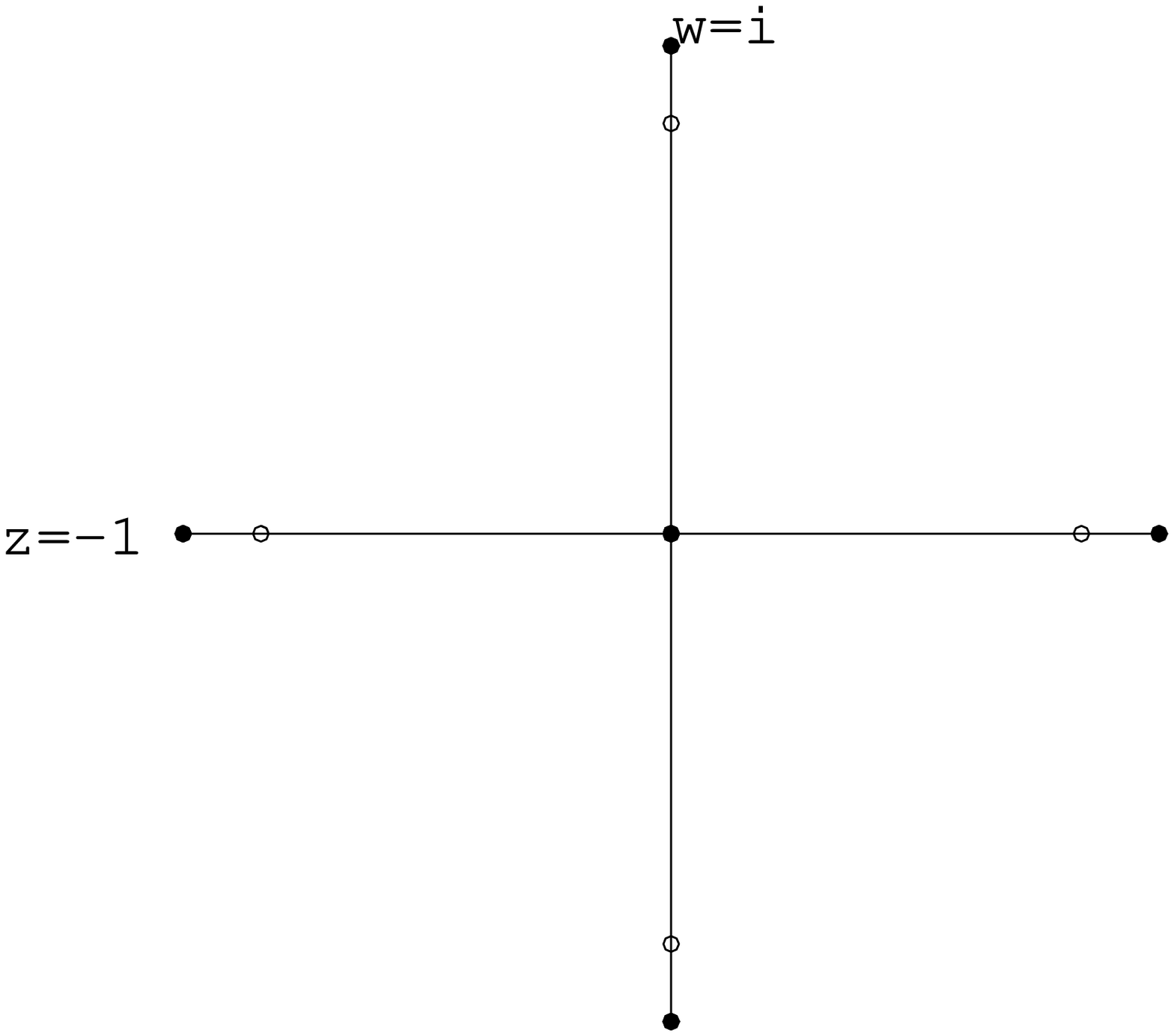}
}

\begin{center}
{\sf Figure 2\\ The normalized dessins of $g$.}
\end{center}
\end{minipage}

\end{center}
\end{figure}

We now prove this combinatorially by appealing to Theorem
\ref{thm:monodromy-tower} with $n=2$, i.e. we shall show that
$f^{\circ 2}$ and $g^{\circ 2}$ have nonisomorphic monodromy groups.
To ease the computation we exploit the fact that $f^{\circ 2} = q
\circ F$ and $g^{\circ 2} = q \circ G$.  Since $q$ is defined over
$\Q$, it suffices to show that $F$ and $G$ have nonisomorphic
monodromy groups.  The (non-clean) dessins of the degree 32
polynomials $F$ and $G$ are shown in Figures 3 and 4.  An obvious
difference is the presence of dihedral symmetry in the dessins of $F$
which is absent in that of $G$.  Note that as abstract one-complexes,
the dessins of $F$ and $G$ are homeomorphic; they differ only in the
way in which they are embedded in the plane.  The dessins of $f$ and
$g$ are obtained from those of $F$ and $G$ by replacing each white
vertex with a black vertex, and then bisecting each edge with a white
vertex.  Labelling the edges of ${\bf D}_F, {\bf D}_G$ more or less
arbitrarily and writing down the elements of $S_{32}$ corresponding to
the generators $\sigma^0, \sigma^1$ for each map, a brief (1-second)
computation in Maple shows that the order of $\Mon(F)$ is $2^{16}$
while the order of $\Mon(G)$ is $2^{18}$, proving our claim. 
Alternatively, one can work directly with ${\bf D}_{f^{\circ 2}}, {\bf
D}_{g^{\circ 2}}$ and determine that the orders of their monodromy
groups are respectively $2^{27}$ and $2^{29}$.  

\begin{figure}[htp]
\begin{minipage}{.45\hsize}
\label{fig:dessinF.ps}
\centerline{%
\epsfysize=.95\hsize
\epsfbox{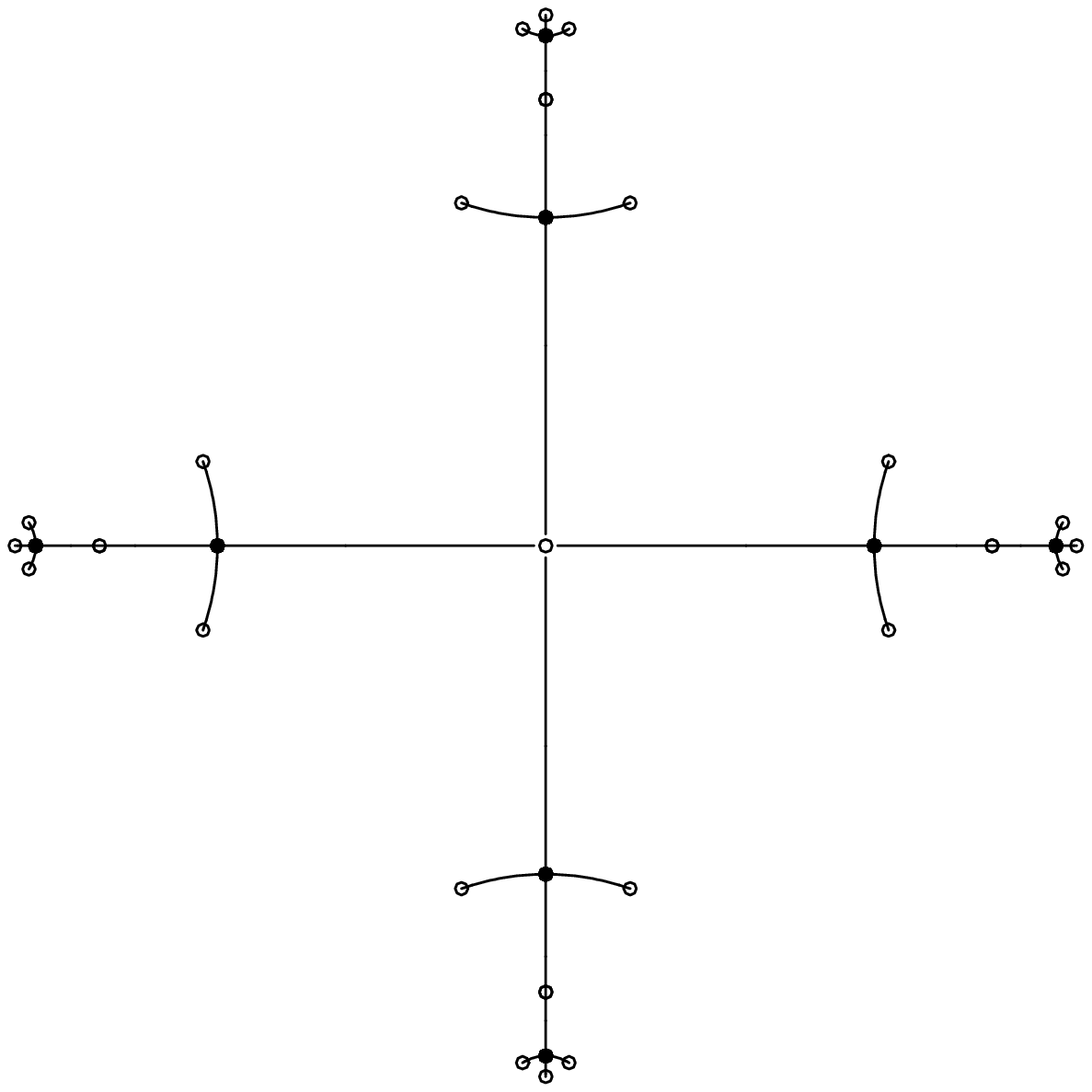}
}
\begin{center}
{\sf Figure 3\\ The dessins of $F$.}
\end{center}
\end{minipage}
\hfill
\begin{minipage}{.45\hsize}

\label{fig:dessinG.ps}
\centerline{%
\epsfysize=\hsize  
\epsfbox{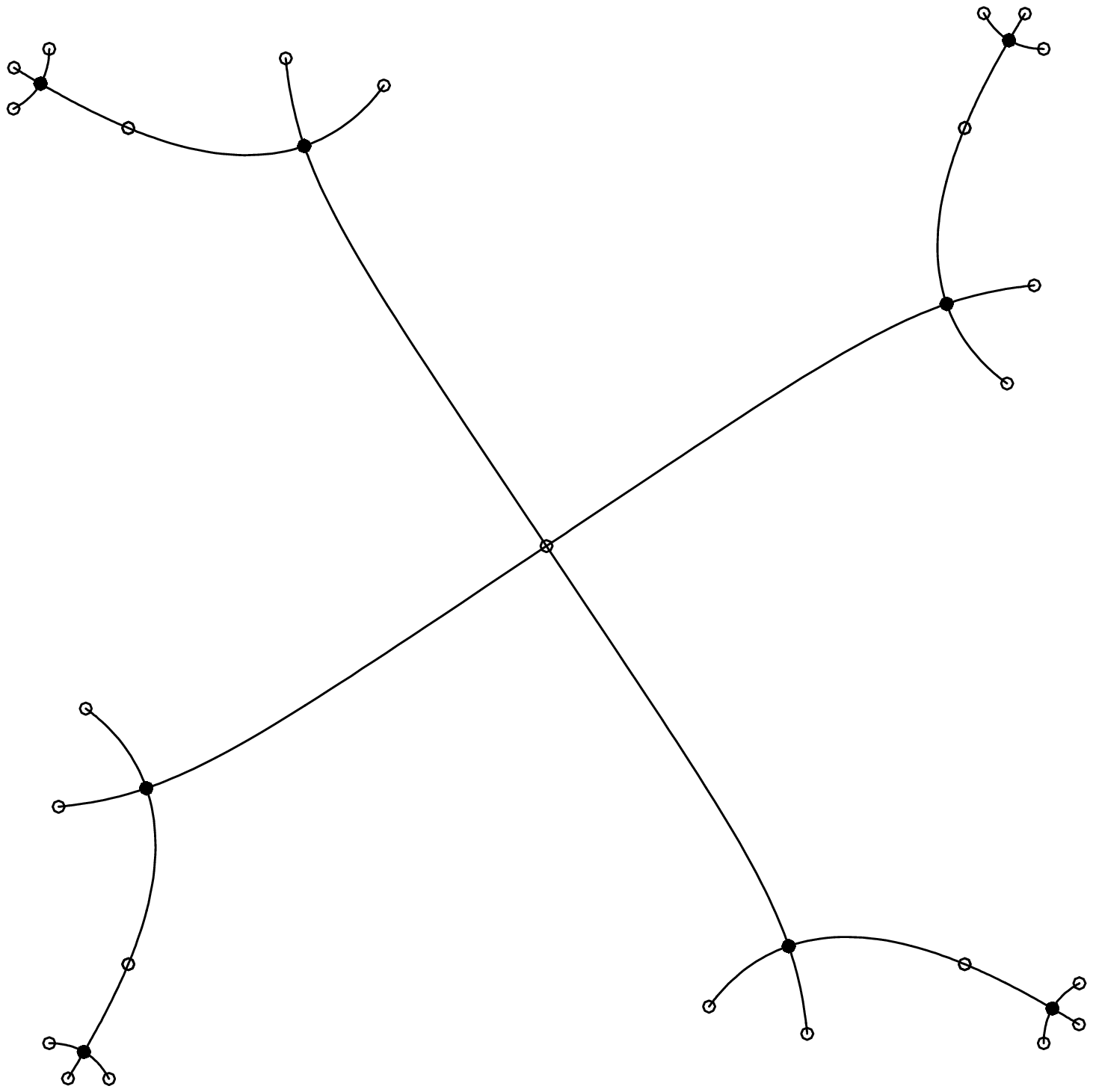}
}
\begin{center}
{\sf Figure 4\\ The dessins of $G$.}
\end{center}
\end{minipage}

\begin{minipage}{.45\hsize}
\centerline{%
\epsfysize=\hsize 
\epsfbox{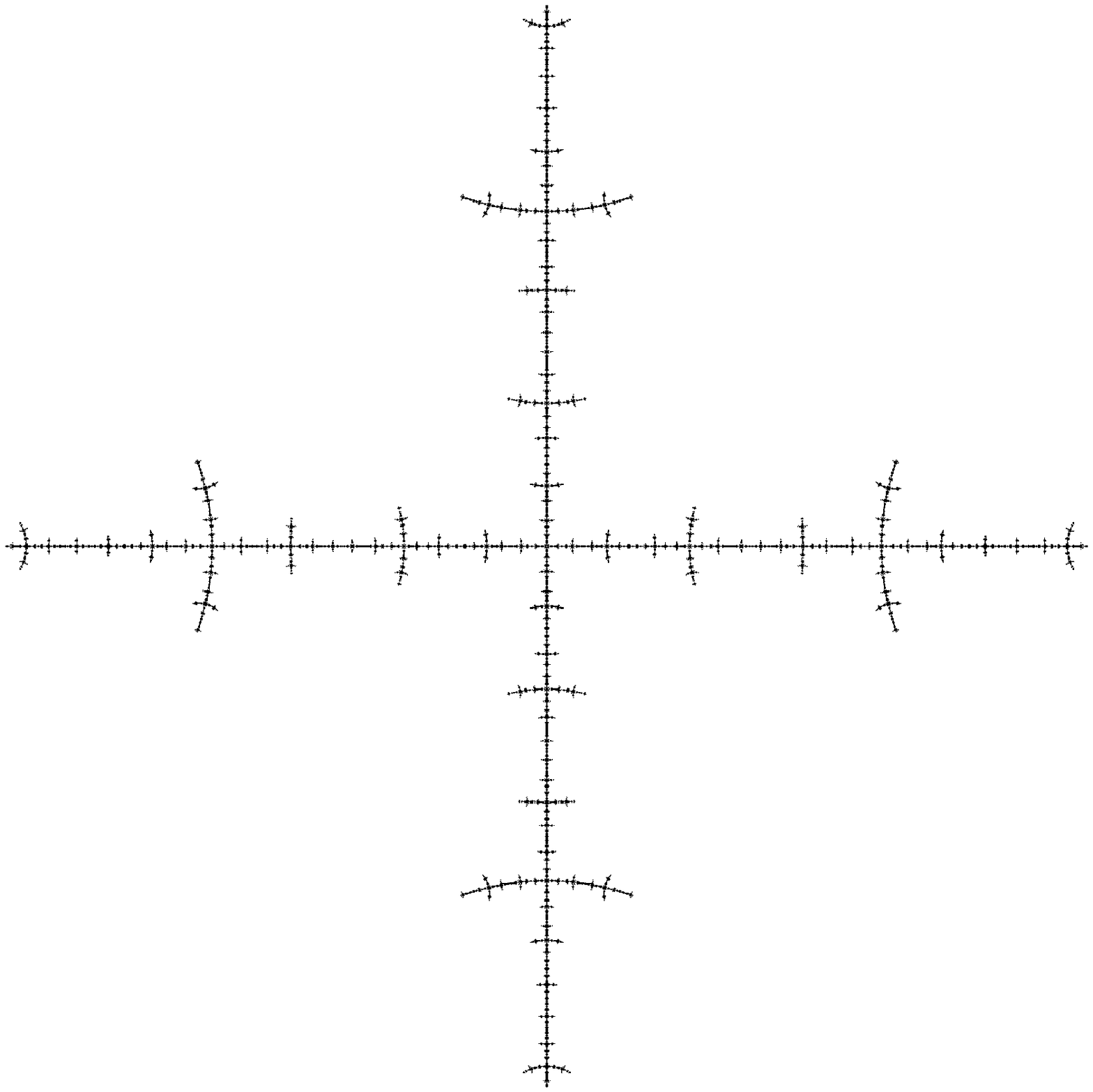}
}
\begin{center}
{\sf Figure 5\\ The Julia set of $f$.}
\end{center}
\end{minipage}
\hfill
\begin{minipage}{.45\hsize}

\centerline{%
\epsfysize=\hsize 
\epsfbox{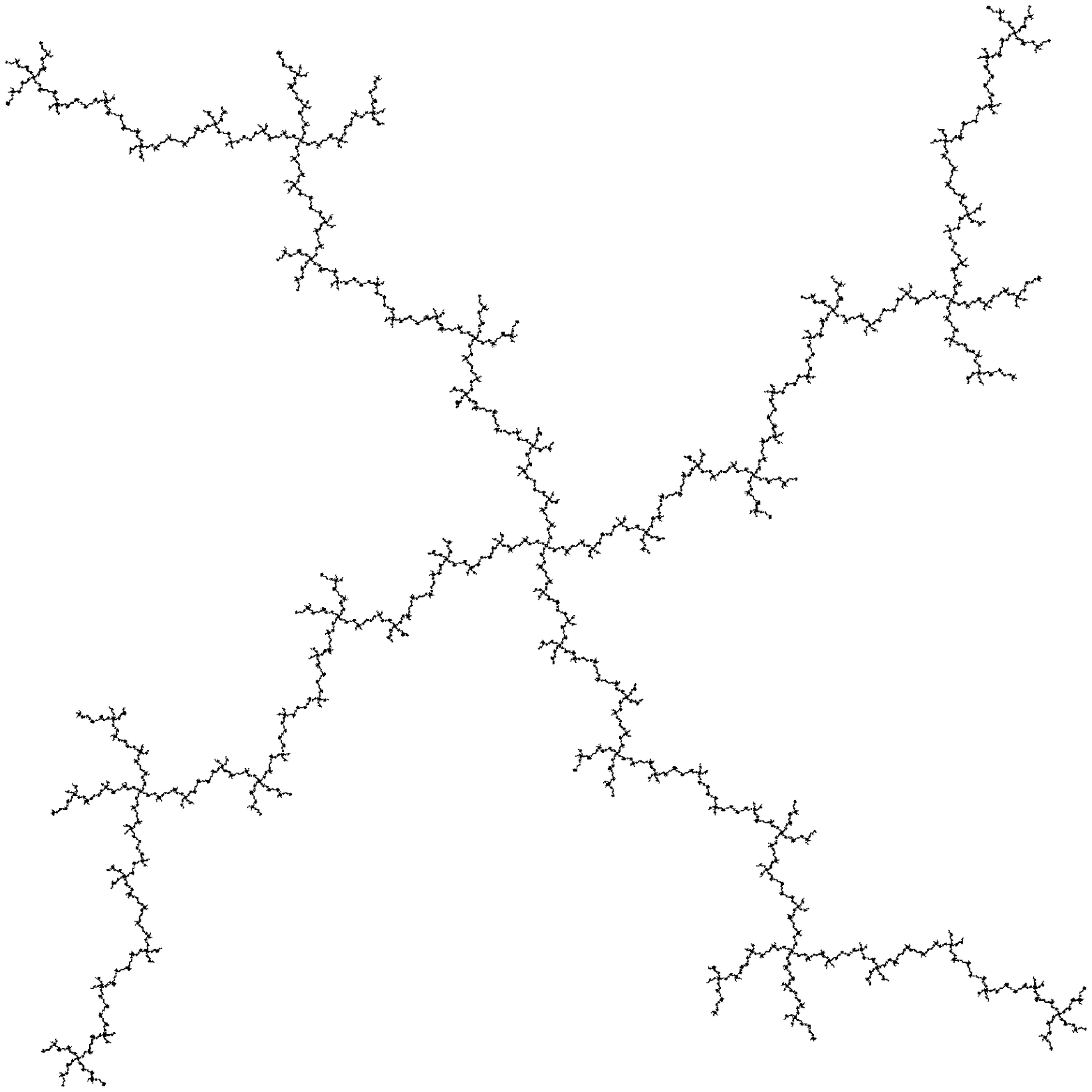}
}
\begin{center}
{\sf Figure 6\\ The Julia set of $g$.}
\end{center}
\end{minipage}
\end{figure}

\noindent{\bf Remark:} As dynamical systems, elements of $XDBP$ are
highly expanding with respect to a suitable (orbifold) metric on a
neighborhood of their Julia sets.  This, and standard arguments from
complex dynamics, can be used to prove that in fact the dessins of
$f^{\circ n}$ converge exponentially fast to the Julia set of $f$.  The
Julia sets of $f$ and $g$ are depicted in Figures 5 and 6.
\gap

\noindent{\bf Recursive formulae for monodromy generators.}  Let $f
\in XDBP$ and let $\sigma_0 ^1 :=\sigma_0, \sigma_1 ^1 :=\sigma_1$ be
the generators for $\Mon(f)$ defined by the action of simple
counterclockwise-oriented loops
$\alpha_0, \alpha_1 \in \pi_1 (\C - \{0, 1\}, b) $ on the fiber
$f^{-1}(b)$ (equivalently, on the set of edges of ${\bf D}_f$).  Here,
we derive recursive formulae for the generators $\sigma_0 ^n, \sigma_1
^n$ of $\Mon(f^{\circ n})$.

Since $f$ is extra-clean, each of the vertices $0$ and $1$ is incident
to exactly one edge of ${\bf D}_f$.  Let $E = \{\epsilon_0,
\epsilon_1, ...\}$ denote the set of edges of ${\bf D}_f$, where
$\epsilon_0$ is the unique edge incident to $0$ and $\epsilon_1$ is
the unique edge incident to $1$; see Figure 7.  We will show that
there is a canonical identification of the set $E^n$ of edges of ${\bf
D}_{f^{\circ n}}$ with the $n$-fold Cartesian product $E \times E
\times ... \times E$ such that the following theorem holds:

\begin{figure}
\centerline{%
\epsfysize=.4\hsize 
\epsfbox{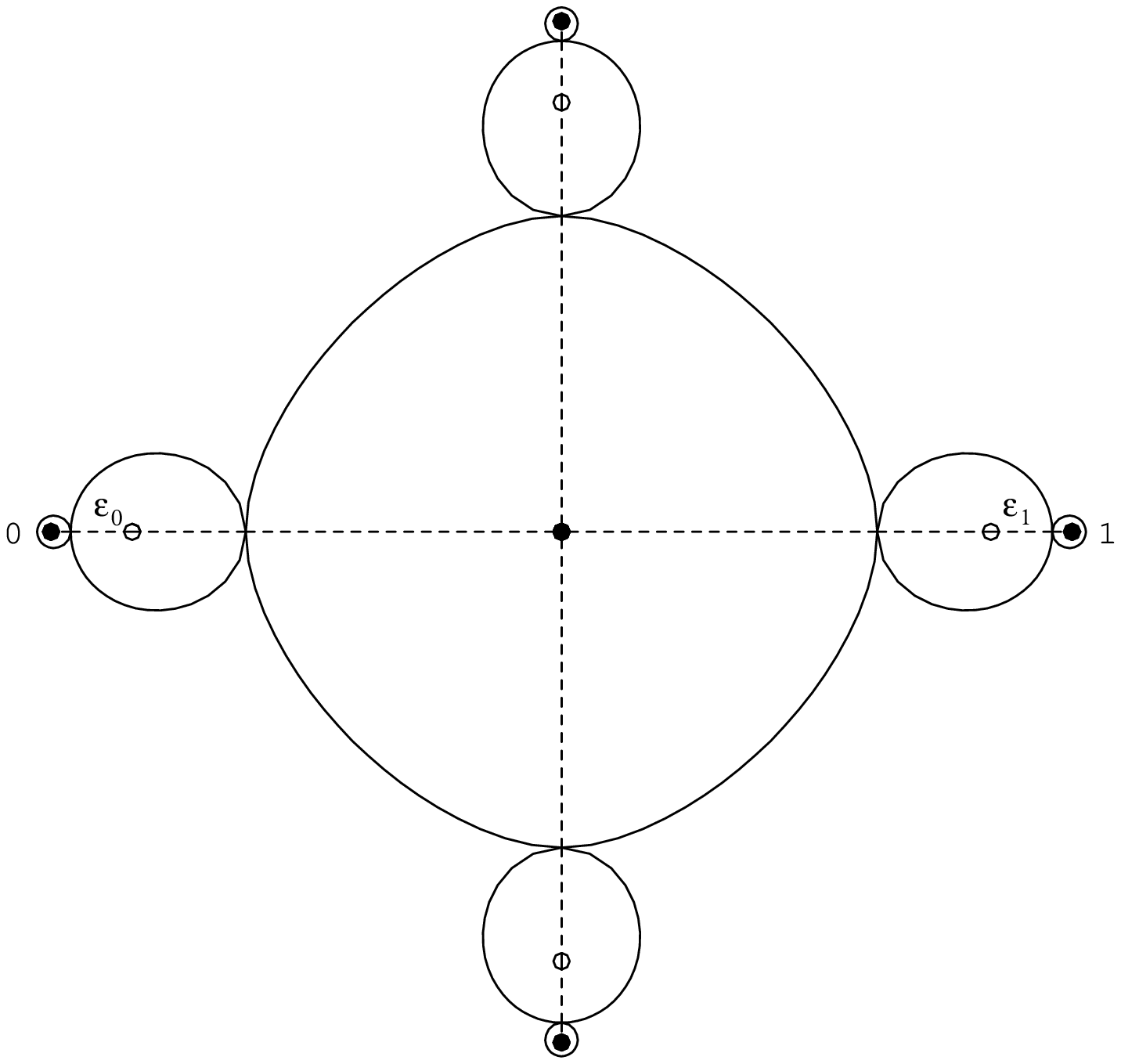}
}
\begin{center}
{\sf Figure 7}
\end{center}
{\sf Here, the lifts of $\alpha_0$ and $\alpha_1$ are the
loops going once around the black and white vertices, respectively.
The edges $\epsilon_0, \epsilon_1$ are indicated.  Note that 
the lift $\alpha_0 ^1 (\epsilon_1)$ of $\alpha_0$ based at
$\epsilon_1$ goes once around $1$, and the lift $\alpha_0 ^1
(\epsilon_0)$ of $\alpha_0$ based at $\epsilon_0$ goes once around
$0$.  The dashed tree is the dessins of the map $f$ of Figure 5.  We have
$U^0 = \C - (-\infty, 0] \union [1, +\infty)$, and all lifts of
$\alpha_0, \alpha_1$ except $\alpha_0^1 (\epsilon_0), \alpha_0 ^1
(\epsilon_1)$ are contained in $U^0$.  
}
\end{figure} 

\begin{thm}
\label{thm:recursive}
For all $n \geq 2$, 
\[ 
\sigma^n _1 (\epsilon_{i_1}, \epsilon_{i_2}, ..., \epsilon_{i_n}) = (\epsilon_{i_1}, \epsilon_{i_2},
..., \epsilon_{i_{n-1}}, \sigma_1 (\epsilon_{i_n})),
\]
and 
\[ \sigma^n _0 (\epsilon_{i_1}, \epsilon_{i_2}, ..., \epsilon_{i_n}) = 
\left\{ 
\begin{array}{lcl}
\;\;\;\;\;\;\;\;\;\;(\epsilon_{i_1}, \epsilon_{i_2},..., \epsilon_{i_{n-1}}, \sigma_0 (\epsilon_{i_n})) & \mbox{if}& \epsilon_{i_n}
\neq \epsilon_0, \epsilon_1 \\
(\sigma_{1}^{n-1}(\epsilon_{i_1}, \epsilon_{i_2}, ...,\epsilon_{i_{n-1}} ),
\epsilon_{i_n}) & \mbox{if}&
\epsilon_{i_n} = \epsilon_1 \\
(\sigma^{n-1} _0 (\epsilon_{i_1}, \epsilon_{i_2}, ..., \epsilon_{i_{n-1}}), \epsilon_{i_n}) & \mbox{if}&
\epsilon_{i_n} = \epsilon_0
\end{array}\right. 
\]
\end{thm}

\noindent We will make use of a convenient Markov partition for the
dynamics of $f$.  \gap

\noindent{\bf The partitions $U^n$.}  Results from holomorphic
dynamics imply the existence of two canonical disjoint arcs $\gamma_0,
\gamma_1$ joining $0$ and $1$ respectively to infinity, called {\it
external rays}, such that $f(\gamma_0) = \gamma_0$ and $f(\gamma_1) =
\gamma_0$.  Set $\Gamma^0 = \gamma_0 \union \gamma_1$ and let $U^0 =
\C - \Gamma^0$; note that $U^0$ is simply-connected and contains no
critical values of $f^{\circ n}$.  Since, conceivably, $\Gamma^0$ may
intersect $(0,1)$, we choose any arc $\lambda$ joining $0$ and $1$,
avoiding $\Gamma^0$, and passing through a rational basepoint $b \in
(0,1)$ for the construction of the dessins ${\bf D}_{f ^{\circ n}}$.
For $n \geq 2$, define inductively $U^n := f^{-1}(U^{n-1})$.  Then
each connected component $U$ of $U^n$ is an open disc mapping
biholomorphically onto its image (which is a connected component of
$U^{n-1}$), and is also contained in a unique connected component of
$U^{n-1}$.  Clearly, the set of connected components of $U^n$ can be
canonically identified with either the set $E^n$, or with the full
preimage $f^{-n}(b)$.  Moreover, since $U^1 \subset U^0$, and each of
the $d$ connected components of $U^1$ maps biholomorphically to $U^0$,
each of the sets $U^n, E^n, f^{-n}(b)$ have $d^n$ elements, where
$d=\deg(f)$.  \gap

\noindent{\bf Itineraries.}  We now show that $E^n$ may be canonically identified with the $n$-fold Cartesian product $E \times E \times ... \times E$.  Given a connected component $U$ of $U^n$ and a subset $X \subset U$, we may form its {\it itinerary} $\iota(X)$ with respect to the partition $U^1$ by setting 
\[ \iota(X) = (\epsilon_{i_1}, \epsilon_{i_2}, ..., \epsilon_{i_n} ) \]
where 
\[ \epsilon_{i_k} = \mbox{the unique component of $U^1$ containing 
$f^{\circ k}(X)$ }, \;\; k=0,1,..., n-1.\] Applying this construction
with $X$ equal to an edge of ${\bf D}_{f^{\circ n}}$ (equivalently,
with an element of $f^{-n}(b)$) we have that $\iota$ determines a
bijection 
\[ \iota: E^n \to \underbrace{E \times E \times ... \times E}_{n}.\]  
Since this identification is canonical, we will use the notation $E^n$
to denote the $n$-fold Cartesian product of $E$ with itself.  
\gap

\noindent{\bf Canonical associated maps.}  The map $f: U^n \to U^{n-1}$ and the
inclusion map $j^n:~U^n~\hookrightarrow~U^{n-1}$ induce maps $f: E^n
\to E^{n-1}$ and $j^n : E^n \to E^{n-1}$ 
given by
\[ f((\epsilon_{i_1}, \epsilon_{i_2}, ..., \epsilon_{i_{n-1}},
\epsilon_{i_n})) =
     (\epsilon_{i_2}, ..., \epsilon_{i_{n-1}}, \epsilon_{i_n}) \]
and
\[ j^n ((\epsilon_{i_1}, \epsilon_{i_2}, ..., \epsilon_{i_{n-1}},
\epsilon_{i_n}) = 
     (\epsilon_{i_1}, \epsilon_{i_2}, ..., \epsilon_{i_{n-1}}) \]
i.e. are the left- and right-shift maps, respectively.
\gap

\noindent{\bf Proof of Theorem.}  To ease notation, set 
\[ (\epsilon) := (\epsilon_{i_1}, \epsilon_{i_2}, ...,
\epsilon_{i_{n-1}},
\epsilon_{i_n}). \] Given $(\epsilon)$, thought of as an element of
$f^{-n}(b)$, we denote by $\alpha_0 ^n ((\epsilon)), \alpha_1 ^n
((\epsilon))$ the lifts of
$\alpha_0, \alpha_1$ under $f^{\circ n}$ based at $((\epsilon))$.
Note that for $1 \leq k \leq n-1$,
\[ f^{\circ k}(\alpha_0 ^n ((\epsilon))) = \alpha_0 ^{n-k}(\epsilon_{i_{k+1}}, ..., \epsilon_{i_n}) \]
and 
\[ f^{\circ k}(\alpha_1 ^n ((\epsilon))) = \alpha_1 ^{n-k}(\epsilon_{i_{k+1}}, ..., \epsilon_{i_n}). \]
\gap

\noindent{\bf Case of $\sigma_1 ^n$.}  We have  
\begin{equation}
\label{eqn:c}
f^{\circ (n-1)}(\alpha_1 ^n ((\epsilon)) = \alpha_1 ^1 (\epsilon_{i_n}) \subset U^0
\end{equation}
where the inclusion at right follows since $f$ is clean; cf. Figure 7.
This implies that
\begin{equation}
\label{eqn:a}
\alpha_1 ^n ((\epsilon)) \mbox{ is contained in a unique component of
$U^{n-1}$.}  
\end{equation} 
By definition, $\sigma_1 ^n ((\epsilon))$ is the endpoint of the curve
$\alpha_1 ^n ((\epsilon))$ and $((\epsilon))$ is its starting point.
Thus by (\ref{eqn:a}), we have
\begin{equation}
\label{eqn:b}
j^n (\sigma_1 ^n ((\epsilon)) = j^n ((\epsilon)).
\end{equation}
But (\ref{eqn:c}) also implies that the image of the endpoint of
$\alpha_1 ^n ((\epsilon))$ under $f^{\circ (n-1)}$ is the endpoint of
$\alpha_1 ^1 (\epsilon_{i_n})$, which is $\sigma_1 ^1
(\epsilon_{i_n}) = \sigma_1 (\epsilon_{i_n})$ by the definition of
$\sigma_1$. This, and (\ref{eqn:b}), prove the formula.  \gap

\noindent{\bf Case of $\sigma_0 ^n$.}  The formula for the case when
$\epsilon_{i_n} \neq \epsilon_0, \epsilon_1$ follows from the same
argument as above, using the fact that since $f$ is extra-clean,
$\alpha_0 ^1 (\epsilon_{i_n}) \subset U^0$ if $\epsilon_{i_n} \neq
\epsilon_0, \epsilon_1$.

Now suppose that $\epsilon_{i_n} = \epsilon_1$.  Since $f$ is
extra-clean, $1$ maps to $0$ under $f$ by local degree one.  Let us
identify $\epsilon_1$ with the unique preimage of $b$ which lies on
the edge incident to $1$.  Then the loop $\alpha_0$ lifts (in the
sense of covering spaces) under $f$ to a simple loop $\alpha_1 ^1
(\epsilon_1)$ which goes exactly once counterclockwise around the
point $1$; see Figure 7.  It then follows that $\sigma_0 ^n $ maps
$f^{-(n-1)}(\epsilon_1)$ to itself bijectively.  Let $\omega$ be an
oriented embedded arc in $U^0$ going from $\epsilon_1$ to $b$.
Lifting $\omega$ determines a bijection
\[ \tau^{n-1}: f^{-(n-1)}(\epsilon_1)  \to f^{-(n-1)}(b)\]
and any two such arcs $\omega, \omega'$ determine the same bijection
since $U^0$ is simply-connected.  The loops
$\overline{\omega}*\alpha_1 * \omega$ (i.e. do $\omega$, then
$\alpha_1$, then the reverse of $\omega$) and $\alpha_0 ^1
(\epsilon_1)$ are both based at $\epsilon_1$ and are homotopic, hence
\[ (\tau^{n-1})^{-1} \circ \sigma_1 ^{n-1} \circ \tau^{n-1} = \sigma_0
^n | f^{-(n-1)}(\epsilon_1). \] 
Note that since $\omega \subset U^0$, any lift of $\omega$ under
$f^{\circ (n-1)}$ is contained in a unique connected component of
$U^{n-1}$.  It follows that $\epsilon$ and $\tau^{n-1}(\epsilon)$ have
itineraries which agree to the first $n-1$ places.  Since 
\[ f^{-(k-1)}(\epsilon_1) = E^{n-1} \times
 \{\epsilon_1\}\]
and 
\[ f^{-(n-1)}(b) = E^{n-1},\]
the map $\tau^{n-1}$ is thus that given by
\[ (\epsilon_{i_1}, ..., \epsilon_{i_{n-1}}, \epsilon_1) \mapsto
(\epsilon_{i_1}, ..., \epsilon_{i_{n-1}}).\]  
That is, 
\[ \tau^{n-1} = j^{n-1} | E^{n-1} \times \{\epsilon_1\}.\]
From this, the result follows immediately.

Finally, the case when $\epsilon_{i_n} = \epsilon_0$ is exactly
analogous--since $f$ is extra-clean, $0$ maps to $0$ by local degree
one, hence the loop $\alpha_0$ lifts to a loop $\alpha_0 ^1
(\epsilon_0)$ going exactly once counterclockwise around 0; see Figure 7.
\qed

\noindent{\bf Remarks.}  The dessins ${\bf D}_{f^{\circ n}}$ may be
inductively constructed as follows.  Recall that by cleanness, each
edge $(\epsilon)^{n-1}$ of ${\bf D}_{f^{\circ (n-1)}}$ is incident to
exactly one black and one white vertex.  To obtain ${\bf D}_{f^{\circ
n}}$, replace each closed edge $(\epsilon)^{n-1}$ of ${\bf
D}_{f^{\circ (n-1)}}$ with a copy of ${\bf D}_f$ by gluing 0 to the
black vertex of $(\epsilon)^{n-1}$ and 1 to the white vertex of
$(\epsilon)^{n-1}$.  This makes synthetic drawing of the dessins of
$f^{\circ n}$ easy.  In contrast, the expanding nature of maps in
$XDBP$ implies that fine detail is rapidly lost in exact drawings, as
Figures 3 and 4 indicate.

\section{Algebraic invariants}

In this section, we formulate some algebraic invariants attached to
a dynamical Belyi polynomial $f$.  We wish these invariants to be
dynamically meaningful, i.e. they should be invariants of the affine
conjugacy class of $f$.  Our discussion of fields of moduli is
drawn from Silverman \cite{Si}. We begin with some background.
\gap

\noindent{\bf The group} $\Gamma$. The group $\Gamma = \GGamma$
inherits a natural topology known as the {\it Krull topology}, and one
has the following as a special case of {\it Krull's theorem}
\cite{Jac}.  Given a subgroup $G < \Gamma$ we denote by $\Inv(G)$
the subfield of $\Qbar$ fixed by $G$.   

\begin{thm}{\rm\bf (Krull)}  Let $\GGG$ be the set of closed
subgroups $G$ of $\Gamma$ and $\FFF$ the set of subfields of $\Qbar$.
Then the maps 
\[ \GGG \ni G \to \Inv(G) \in \FFF\]
and
\[ \FFF \ni F \to \Gal(\Qbar/F) \in \GGG\]
are inverses and order-reversing.  A subgroup $G \in \GGG$ is normal
in $\Gamma$ if and only if $F = \Inv(G)$ is Galois over $\Qbar$, and
if this is the case, then $\Gal(F/\Q) \cong \Gamma/\Gal(\Qbar/F) $.  
\end{thm}
\gap

In what follows, we let $\Rat_d$ denote the set of all
rational maps of a given degree $d \geq 2$ with coefficients in
$\Qbar$.  We remark that $\Gamma$ acts on $\Rat_d$ by
twisting coefficients, as with the case of polynomials.  
\gap

\noindent{\bf Fields of definition.}  Given $f \in \Rat_d$ and a field 
extension $K$ of $\Q$, we say that $f$ is {\it defined over
$K$} if its coefficients are in $K$, and that $K$ is a {\it field of
definition} for $f$.  It can easily be shown that $f$ is defined over
$K$ if and only if $f^{\sigma} = f$ for all $\sigma \in \Gal(\Qbar/K)$.

\begin{defn}
Given $f \in \Rat_d$, the {\bf\rm coefficient field
$K_{\Coeff} (f)$ } of $f$ is the field generated by its coefficients.
\end{defn}
\gap

\noindent{\bf Fields of moduli.}  If $f, g \in \Rat_d$ are
affine conjugate, i.e. $f = A g A^{-1}, A \in \PSL_2 (\C)$, 
then necessarily $A \in \PSL_2 (\Qbar)$ since
e.g. $A$ must send periodic points to periodic points:  there are
always at least three such points, which are necessarily algebraic,
and $A$ is determined by where it sends three points.  Hence  
\[ f^{\sigma} = A^{\sigma} g^{\sigma} (A^{\sigma})^{-1}\]
and so the action of $\Gamma$ on $\Rat_d$ descends to an
action on the {\it moduli space} 
\[ M_d = \Rat_d /\mbox{conjugation}.\]
Denoting elements of $M_d$ by $\varphi$ and the action by $\varphi
\mapsto \varphi^{\sigma}$, we note that by definition, for any $f \in
\varphi$, 
\[ \varphi^{\sigma} = \mbox{the conjugacy class of }f^{\sigma}.\]
Given $\varphi \in M_d$, we denote the stabilizer of $\phi$
under the action of $\Gamma$ by 
\[ \Stab_{\Gamma} (\varphi) = \{ \sigma \in \Gamma | \varphi^{\sigma}
= \varphi \}.\]

\begin{defn}
Let $\varphi \in M_d$.  The {\rm field of moduli
$K_{\Moduli}(\varphi)$ of
$\varphi$} is defined by \[ K_{\Moduli}(\varphi) =
\Inv(\Stab_{\Gamma}(\varphi)).\] We say that a field $K \subset \Qbar$ is a
{\bf\rm field of definition for $\varphi$} if $\varphi$ contains an element $f$
defined over $K$.
\end{defn}

The field of moduli of $\varphi \in M_d$ is contained in any field of
definition of $\phi$.  To see this, choose any $f \in \varphi$ defined
over $K \subset \Qbar$.  Then 
$f^{\sigma} = f$ for all $\sigma \in \Gal(\Qbar/K) < \Gamma$.  Hence
$\varphi^{\sigma} = \varphi$ for all $\sigma \in \Gal(\Qbar/K)$ and so
$\Gal(\Qbar/K) < \Stab_{\Gamma}(\varphi)$. By Krull's theorem, $K =
\Inv (\Gal(\Qbar/K)) \supset \Inv(\Stab_{\Gamma}(\varphi)) =
K_{\Moduli}(\varphi)$.  

The question of whether the field of moduli of a dynamical system
$\varphi \in M_d$ is also a field of definition is quite subtle.
Silverman \cite{Si} has shown that if $d$ is even, or if $\varphi$
contains a polynomial, then the field of moduli is always a field of
definition.  On the other hand, he shows that e.g.  the conjugacy
class $\phi$ of the map
\[ f(z) = i \left( \frac{z-1}{z+1} \right)^3 \] 
has field of moduli equal to $\Q$ (since $f$ is M\"obius conjugate to
$\overline{f}$), but nonetheless is not definable over $\Q$.
Couveignes \cite{Cou} has studied the analogous question for
non-dynamical Belyi polynomials (i.e. where the equivalence relation
is $f\sim g$ if $f=g\circ A$, $A \in \Aut(\C)$), and gives an example
of a polynomial whose field of moduli is not a field of definition.

Recall that $DBP$ is a set of polynomials, not ``modded out'' by an
equivalence relation.  In our case, it turns out that the coefficient
field of an element $f \in DBP$ is an intrinsic quantity of the
conjugacy class of $f$:

\begin{thm}
\label{thm:moduli}
Let $\varphi \in M_d $ be the conjugacy class of the degree $d$ map $f
\in DBP$.  Then $K_{\Coeff}(f) = K_{\Moduli}(\varphi)$.
\end{thm}

\pf The stabilizer in $\Gamma$ of $f \in DBP$ is a closed subgroup of $\Gamma$.
Since no two elements of $DBP$ are affine conjugate, the stabilizers of $f$
under the action of $\Gamma$ on $DBP$ and of $\varphi$ under the action of
$\Gamma$ on $M_d$ coincide.  Hence their invariant subfields of $\Qbar$ are the
same, and so 
\[ K_{\Moduli}(\varphi) = \Inv(\Stab_{\Gamma}(\varphi)) =
\Inv(\Stab_{\Gamma}(f)).\]
Since the field of moduli is always contained in any field of
definition, 
\[ K_{\Moduli}(\varphi) \subset K_{\Coeff}(f).\]
To prove the other direction, it suffices to show that $\Inv(\Stab_{\Gamma}(f))
\supset K_{\Coeff}(f)$, and since $\Stab_{\Gamma}(f)$ is closed this
is equivalent  to showing that $\Stab_{\Gamma}(f) <
\Gal(\Qbar/K_{\Coeff}(f))$, by Krull's theorem.  But this is clear, since
$\sigma \in \Gamma$ is the identity on $K_{\Coeff}(f)$ if and only if it is
the identity on the generating set of $K_{\Coeff}(f)$, i.e. on the
coefficients of $f$.  
\qed

\begin{cor}If $f \in DBP$ and $\varphi$ is the conjugacy class of $f$,
then $\varphi$ is definable over $\Q$ if and only if $f \in \Q[z]$.
\end{cor}

Since elements of $\DBP$ admit no nontrivial automorphisms, the set
$Z_f$ of zeros, the set $O_f$ of ones, and the set $C_f \supset O_f$
\footnote{This condition is guaranteed by the assumption of
cleanness.}  of critical points are all dynamically distinguished
subsets of $\Qbar \subset \C$.  Hence they form invariants of the
conjugacy class of $f$.  For a subset $S$ of $\Qbar$ and a field $K$,
we denote by $K(S)$ the extension of $K$ generated by $K$ and $S$.
The above fields are related as follows:

\begin{thm}  
Let $f \in DBP$ and set $K=K_{\Coeff}(f)$.  Then we have the
following collection of field extensions, each of which is Galois: 
\end{thm}

\begin{center}
\begin{picture}(160,160)
\put(80,0){\makebox(0,0){$K$}}
\put(40,40){\makebox(0,0){$\Q(O_f)$}}
\put(120,40){\makebox(0,0){$\Q(Z_f)$}}
\put(40,80){\makebox(0,0){$\Q(C_f)$}}
\put(80,120){\makebox(0,0){$\Q(Z_f \union O_f)$}}
\put(40,35){\line(1,-1){30}}
\put(120,35){\line(-1,-1){30}}
\put(40,75){\line(0,-1){30}}
\put(80,115){\line(-1,-1){30}}
\put(80,115){\line(1,-2){35}}
\end{picture}
\gap

{\sf Figure 8}

\end{center}

We first establish 

\begin{lemma}
Let $f \in DBP$ and let $K = K_{\Coeff}(f)$. Then $K \subset \Q(Z_f) \cap \Q(O_f) $.  
\end{lemma}

\pf (of Lemma). We have 
\[ f(z) \;\; = \;\; c \cdot \prod_{z_i \in Z_f} (z-z_i)^{d_i} \;\; =\;\; c \cdot
\prod_{w_j \in O_f} (z-w_j)^2 \] 
for some $c \in \Qbar$.  It suffices to show $c
\in \Q(Z_f)$ and $c \in \Q(O_f)$.  Suppose $c \not\in \Q(Z_f)$.  Let $E$ be a
Galois extension of $\Q(Z_f)$ containing $c$ and let $\sigma \in
\Gal(E/\Q(Z_f))$ satisfy $\sigma(c) \neq c$.  Then $f \neq f^{\sigma}$ but $Z_f
= Z_{f^{\sigma}}$, violating Theorem \ref{thm:zofix}.  The case of $O_f$ is
similar.
\qed

\pf (of Theorem). By the Lemma, the diagram in the Theorem is equivalent to 

\begin{center}
\begin{picture}(160,160)
\put(80,0){\makebox(0,0){$K$}}
\put(40,40){\makebox(0,0){$K(O_f)$}}
\put(120,40){\makebox(0,0){$K(Z_f)$}}
\put(40,80){\makebox(0,0){$K(C_f)$}}
\put(80,120){\makebox(0,0){$K(Z_f \union O_f)$}}
\put(40,35){\line(1,-1){30}}
\put(120,35){\line(-1,-1){30}}
\put(40,75){\line(0,-1){30}}
\put(80,115){\line(-1,-1){30}}
\put(80,115){\line(1,-2){35}}
\end{picture}
\gap

{\sf Figure 9}

\end{center}

As extensions of $K$, the fields $K(O_f)$, $K(Z_f)$, $K(C_f)$, and
$K(Z_f \union O_f)$ are the splitting fields of the separable
polynomials in $K[z]$ given by $f-1$, $f$, $f'$, and $f(f-1)$,
respectively.  Hence all indicated and implied extensions are Galois.
\qed

\noindent{\bf Remark and Example:}  In general, it is possible that $\Q(C_f)$ is
a proper subfield of $Q(Z_f \union O_f)$, as the following example shows. 
Consider the polynomial $f = q \circ g$ where $q(z) = 4z(1-z)$ and $g(z) =
-3z^3(z-4/3)$.  It is easily verified that $f \in DBP$, and that 
\[ \Q(O_f) =\Q(C_f) = \Q(\mbox{Roots of }\; 6z^4 - 8z^3 + 1) \] 
and 
\[ \Q(Z_f) =\Q(\sqrt{2}i).\] 
The Galois group of the above quartic is $S_4$ and its
discriminant is $\Delta = -2^{11}\cdot 3^3$.  Since $S_4$ contains a unique
normal subgroup of index 2, $\Q(C_f)$ contains a unique Galois subfield of
dimension two over $\Q$, namely $Q(\sqrt{\Delta}) = \Q(\sqrt{6}i) \neq
\Q(\sqrt{2}i)$.  Hence $\Q(\sqrt{2}i)$ is not a subfield of $\Q(C_f)$.

\section{Hubbard trees}

In this section, we introduce Hubbard trees and conclude the proof of
our main result: that $\Gamma$ acts faithfully on Belyi-type Hubbard
trees (Theorem \ref{thm:main2}).  We also explicitly relate Hubbard
trees and normalized dessins.  Our proof relies on Thurston's
combinatorial characterization of a certain class of polynomials
regarded as holomorphic dynamical systems, and on work of Poirier
connecting Thurston's characterization with Hubbard trees.  \gap

\noindent{\bf Thurston equivalence of branched coverings.} We relax
the structure on elements of $DBP$ and consider instead the more
flexible setting of {\it continuous} orientation-preserving, finite
degree, postcritically finite branched coverings $F: \C \to \C$.
Following Thurston, two such maps $F, G$ are called {\it
combinatorially equivalent} if there are orientation-preserving
homeomorphisms $\psi_0, \psi_1: \C \to \C$ for which $\psi_0 \circ F =
G \circ \psi_1$ and for which $\psi_0$ is isotopic to $\psi_1$ through
homeomorphisms fixing $P_F$, the postcritical set of $F$.  

\begin{thm} {\rm\bf (Thurston Rigidity)}  If $F: \C \to \C$ is a postcritically
finite branched covering, then $F$ is combinatorially
equivalent to at most one (up to affine conjugacy) polynomial $f$. 
\footnote{More generally, this holds for postcritically finite
branched coverings $F: S^2 \to S^2$, apart from the non-polynomial
Latt\`es examples.}
\end{thm}

As a special case of Thurston's characterization of postcritically
finite rational functions as branched coverings of the sphere to
itself \cite{DH} we have

\begin{thm} 
If $f: \C \to \C$  is a postcritically finite branched covering 
for which $\# P_F = 2$, then $F$ is combinatorially
equivalent to exactly one (up to affine conjugacy) polynomial $f$.
\end{thm}

\noindent{\bf Abstract minimal Hubbard trees.}  A Hubbard tree can be
thought of as giving an almost normal form for a postcritically finite
polynomial map from the plane to itself.  Roughly, it is a finite
planar tree $T$ together with a map $\tau$ from $T$ to itself sending
vertices to vertices (but not necessarily edges to edges) which is
extendable to a branched covering of the plane to itself, satisfies
some reasonable minimality conditions, and is sufficiently expanding
to be combinatorially equivalent to a polynomial.  Our definitions are
slightly condensed versions of those of Poirier \cite{Poi}.  
\gap

\noindent{\bf Definition.}  A {\bf Hubbard tree} \footnote{Poirier
would call these {\it dynamical, homogeneous, expanding, minimal
Hubbard trees.}} is a triple $\bf T = (T, \tau, \delta)$ where 

\begin{enumerate}

\item {\rm \bf Planar tree.} $T$ is a planar tree with vertex set $V$;

\item {\rm \bf Dynamics.} 
$\tau: T \to T$ is a continuous map sending $V$ to $V$ and which
is injective on the closures of edges;

\item {\rm \bf Local degree.} 
$\delta: T \to \{1,2,3,...\}$ is the {\it  local degree
function};

\item {\rm\bf Uniquely extendable.} For each $x \in T$, $\tau$
extends to an orientation-preserving branched covering from a
neighborhood of $v$ to a neighborhood of $\tau(v)$, and which in local
coordinates with $v = \tau(v) = 0$ is given by $\tau(z)
= z^{\delta (v)}$.  If $\delta(v) > 1$ we call $v$ a {\it critical
point} and denote the set of critical points by $C$.  We define the
{\it postcritical set} $P$ by 
\[ P = \bigcup_{n>0} \tau^{\circ n}(C) \subset V.\]

\item {\rm\bf (Nontrivial.)} The {\rm degree} $d$ of $\bf T$ defined
by  \[ d = 1 + \sum_{v \in V} \delta(v) - 1\]
is at least two;

\item {\rm\bf (Homogeneity.) } $C \subset \tau(T)$ (i.e. critical
points must have preimages), and for each $x \in \tau(T)$,
\[ \sum_{\tau(y)=x} \delta(y) = d;\]
(i.e.  all points which have preimages have exactly $d$ preimages,
counting multiplicities);

\item {\rm\bf (Expansion.) } A vertex $v$ is called a {\it Fatou vertex}
if there are integers $n \geq 0, p>0$ for which $\tau^n (v) =
\tau^{n+p}(v) \in C$ (i.e. it lands on a period $p$  critical point);
otherwise $v$ is called a {\it Julia vertex}, and we require that if
$v, v'$ is any pair of distinct, adjacent Julia vertices, then there is an
integer $n > 0$ for which $\tau(v), \tau(v')$ are nonadjacent.

\item {\rm\bf (Minimality.)} Given two such triples $\bf T, \bf T'$,
we say $T \preceq T'$ if there is a dynamically compatible
orientation-preserving embedding of pairs, i.e. an embedding $\phi:
(T,V) \to (T',V')$ such that $\tau'\circ \phi = \phi \circ \tau$ and 
$\delta = \delta' \circ \phi$.  
\end{enumerate}

Finally, we say two Hubbard trees $\bf T, \bf T'$ are {\it isomorphic}
if $T \preceq T'$ and $T' \preceq T$.  Equivalently, two Hubbard trees
are isomorphic if there is an orientation-preserving homeomorphism of
the plane to itself carrying $(T, V)$ to $(T',V')$, conjugating the
dynamics on the vertices, and preserving the local degree functions.

The minimality criterion is obviously necessary if we think of a
Hubbard tree as a reasonable normal form for encoding a branched
covering from the plane to itself.  For example, consider a Hubbard
tree with an edge which is fixed pointwise.  Add a vertex at the
midpoint, and transport this using the dynamics to obtain another new
homogeneous tree.  Finally, add any number of edges emanating from the
new vertex and transport these new edges using the dynamics to obtain
a new homogeneous tree.  Dynamically, these added edges are
extraneous.

\begin{thm} {\bf\rm (\cite{Poi}, Thms. II.4.7 and II.4.8, p. 28)} 
The set of affine conjugacy classes of postcritically finite
polynomials of degree at least two is in bijective correspondence with
the set of isomorphism classes of Hubbard trees.
\end{thm}

We now outline the construction of this bijection.  A Hubbard tree
$\bf T$ determines a postcritically finite branched covering $F: \C
\to \C$.  Expansion and Thurston's characterization of postcritically
finite rational maps as branched coverings guarantees that $F$ is
combinatorially equivalent to a postcritically finite polynomial,
which is unique up to affine conjugacy by Thurston rigidity.  So one
has a well-defined map from Hubbard trees $\bf T$ to conjugacy classes
of polynomials $f_{\bf T}$.  The injectivity of this map follows from
the minimality criterion--without it, one could e.g. simply add
additional orbits of periodic vertices to obtain a new tree which
yields the same polynomial.  As a consequence, two Hubbard trees, when
extended to branched coverings of the plane, yield combinatorially
equivalent coverings if and only if they are isomorphic.

A description of the inverse of this map requires some notions from
complex dynamics.
\gap

\noindent{\bf Julia sets.}  Let $f: \C \to \C$ be a polynomial of
degree $\geq 2$.  We define 
\begin{itemize}
\item the {\it filled Julia set} $K_f = \{z \;| \;f^{\circ n} (z) \not\to
\infty\}$;
\item the {\it Julia set} $J_f = \bdry K_f$;
\item the {\it Fatou set} $F_f = \C - J_f$.
\end{itemize}

If $f$ is postcritically finite, $K_f$ and $J_f$ are connected and
locally connected.  The set $K_f$ has interior if and only if there
are periodic critical points.  Otherwise, $K_f = J_f$ is a {\it
dendrite}, i.e. compact, connected, with no interior, and whose
complement is connected.  In particular, the Julia sets of maps in
$XDBP$ are always dendrites, and given any two points $x, y \in K$,
there is a unique (up to reparameterization) continuous embedding
$\gamma: [0,1] \to K_f$ joining $x$ to $y$.  In the case when $K_f$
has interior, the connected components of the interior of $K_f$ are
bounded components $\Omega$ of $F_f$.  These are Jordan domains whose
closures meet in at most one point.  Each such bounded component
$\Omega$ contains a unique ``center point'' $p$ which eventually maps
onto a critical point.  The pair $(\Omega, p)$ is conformally
isomorphic to the unit disc $(\Delta, 0)$ and thus carries a canonical
foliation by radial arcs $\{re^{2\pi it} | r<1 \}$ called {\it
internal rays}.  Any two points $x,y \in \cl{\Omega}$ are joined by a
unique arc, called a {\it regulated arc}, which is the union of two
closed internal rays.

The tree ${\bf T} _f$ associated to a postcritically finite polynomial
$f$ may now be described as follows.  The underlying tree $T$ is the
smallest subcontinuum of $K_f$ containing $f^{-1}(P_f)$ and whose
intersection with the closure of any bounded Fatou component is either
a regulated arc or a finite collection of points on the boundary. This
is always a topological tree.  The vertex set $V$ one takes to be
$f^{-1}(P_f)$ together with the (necessarily finite) number of points
$v$ for which $T-\{v\}$ has three or more components.  As $\tau$ one
takes the function $f|T$, and it is known that $f(V) \subset V$ and
$f(T) \subset T$.  Finally, as $\delta$ one takes the local degree
function of associated to $f$.  Poirier then shows that if $\bf T$ is
a Hubbard tree and $f_{\bf T}$ the associated polynomial, then the
tree obtained by applying the above construction to $f_{\bf T}$ is a
Hubbard tree equivalent to $\bf T$.

We now make an explicit connection between Hubbard trees and tree
dessins.  
\gap

\goodbreak
\noindent{\bf Definitions.}  We denote by 
\begin{itemize}
\item $[BH]$ the set of isomorphism classes of Hubbard trees $\bf T$
degree at least three for which the postcritical set $P$ has at most
two points, and for which there exists a (necessarily unique) $v \in
P$ such that $\delta(y)=2$ for all $y \in V$ with $\tau(y) = v$,
i.e. the local degree of $\tau$ near $y$ is exactly equal to two.  We
refer to $[BH]$ as the set of isomorphism classes of {\it clean
Belyi-type Hubbard trees}.

\item $[BBC]$ the set of combinatorial equivalence classes of branched
coverings $F: \C \to \C$ of degree at least three for which $P_F$ has
exactly two points, and for which there exists a (necessarily unique)
$v \in P_F$ such that the local degree of $F$ near $y$ is exactly
equal to two for all $y$ with $F(y) = v$.  We refer to $[BBC]$ as the
set of combinatorial equivalence classes of {\it clean Belyi-type
branched coverings} from the plane to itself.
\end{itemize}

We define extra-clean Belyi type Hubbard trees and extra-clean
Belyi-type branched coverings analogously.  Note that the homogeneity
and extendability conditions guarantee that if ${\bf T} \in BH$ and
$F$ is any extension of $\tau$ to a branched covering of the plane to
itself, then $F$ is a clean Belyi-type branched covering.  
\gap

\noindent{\bf Proof of Theorem 3.8.}  
There are natural maps from $DBP$ to $[BBC]$ and $[BHT]$ which send a
dynamical Belyi polynomial $f$ to its combinatorial class as a
branched covering and the isomorphism class of its Hubbard tree,
respectively.  We then have the following commutative diagram, where
all lines indicate bijections:

\begin{center}
\begin{picture}(150,150)
\put(75,125){\makebox(0,0){$[BP^*]$}}
\put(20,75){\makebox(0,0){$[TD^*]$}}
\put(75,75){\makebox(0,0){$DBP$}}
\put(130,75){\makebox(0,0){$[BBC]$}}
\put(75,25){\makebox(0,0){$[BHT]$}}
\put(50,100){\makebox(0,0){II}}
\put(100,100){\makebox(0,0){I}}
\put(50,50){\makebox(0,0){III}}
\put(100,50){\makebox(0,0){IV}}
\put(20,85){\line(1,1){40}}
\put(130,85){\line(-1,1){40}}
\put(20,65){\vector(1,-1){37}}
\put(130,65){\line(-1,-1){37}}

\put(40,75){\line(1,0){20}}
\put(90,75){\line(1,0){20}}
\put(75,35){\line(0,1){30}}
\put(75,85){\line(0,1){30}}
\end{picture}

{\sf Figure 10} 

\end{center}

Here, the diagonal map at upper right is by definition the composition
of the sides, so triangle (I) is commutative by definition.  Triangle
(II) commutes by extension of the Grothendieck correspondence to
normalized clean tree dessins and normalized clean Belyi polynomials.
Triangle (IV) commutes by restricting Poirier's theorems to the case
of clean dynamical Belyi polynomials, Hubbard trees, and branched
coverings.

By Poirier's theorem, $DBP$ is in bijective correspondence with
$[BHT]$.  So we can use this correspondence to define a faithful
action of $\Gamma$ on $[BHT]$.  This applies as well in the
extra-clean case, and Theorem 3.8 is proved.  \qed \gap

\noindent{\bf Hubbard trees and normalized dessins.}  The preceding
proof is somewhat unsatisfactory, since the map sending $f \in DBP$ to
its Hubbard tree (which is a combinatorial object) is defined via the
topology and dynamics of the Julia set, which is geometric.  Below, we
give an alternative description of the arrowed map in Triangle (III)
of Figure 10.  We will construct, given $f \in DBP$ and its normalized
dessin, an associated Hubbard tree $\bf T$ and an extension of the
dynamics to a branched covering which is combinatorially equivalent to
$f$.  By uniqueness of the Hubbard tree, this tree will be {\it the}
Hubbard tree of $f$.

Given any isomorphism class ${\bf D}^*$ of normalized dessin, let $f$
be the corresponding element of $DBP$, let ${\bf D}_f$ be the dessin
of $f$ with vertices $X_0$ (i.e. ${\bf D}_f$ is the preimage of
$[0,1]$ with its tree structure and bicoloring of vertices).  Then
since $f \in DBP$, $f(\{0,1\}) \subset \{0,1\}$, so $f^{-1}(\{0,1\})
\supset \{0,1\}$.  Let $\alpha$ denote the unique subarc of $D_f$
joining $0$ and $1$ in ${\bf D} _f$.  Since any two arcs joining zero
and one in the plane are ambient isotopic fixing endpoints, there is a
homeomorphism $\psi_1: \C \to \C$ be a homeomorphism such that $\psi_1
(0) = 0, \psi_1 (1) = 1$, $\psi_1 ([0,1]) = \alpha$, and $\psi_1$ is
isotopic to the identity through maps fixing $0$ and $1$.  Let $F: \C
\to \C$ be the branched covering $f \circ\psi_1$.  Then $F$ is
combinatorially equivalent to $f$ (take $\psi_0$ to be the identity in
the definition) and indeed this is true for any choice of $\psi_1$.

We now set $\bf T$ to be the Hubbard tree whose underlying tree $T$ is $\psi_1 ^{-1}(D_f)$, whose vertex set $V$ is
$\psi_1 ^{-1}(X_0 (f))$, whose local degree function is the restriction
of the local degree of $F$, and whose dynamics $\tau$ is the
restriction of $F$ to $T$.

We next claim that $\bf T$ is indeed a Hubbard tree.  By construction, edges map
injectively on their interiors and the dynamics is extendable.  Homogeneity is
also clear, since the underlying tree is the full preimage of a set containing
all the finite postcritical set.  Minimality holds since each vertex is the
preimage of a critical value, these are required by the homogeneity condition,
and hence the tree cannot be modified by removing vertices to obtain a smaller
tree with respect to the partial order $\preceq$ on trees.  We now verify the
expansion condition.  Let $v, v'$ be any pair of adjacent Julia vertices, and
suppose that for all $i$, $\tau^{\circ i}(v)$ and $\tau^{\circ i}(v')$ are
adjacent. Note that $\tau^{\circ i}(v) \neq \tau^{\circ i}(v')$ by requirement
\#2 in the definition of Hubbard tree.  Since $P_T$ is finite (in fact, has two
elements), by replacing $v$ and $v'$ with $\tau^{\circ i} (v), \tau^{\circ
i}(v')$ we may assume $v, v'$ are periodic.  By the construction of $\bf T$, we
may assume $v=0$ and $v'=1$.  Then $0$ and $1$ are periodic Julia vertices,
hence neither is a critical point.  That is, both $0$ and $1$ are ends of the
tree ${\bf D}_f$.  Since the degree of $f$ is at least two, any two ends of the
tree ${\bf D}_f$ are nonadjacent, a contradiction.

The above function $[TD^*] \to [BHT]$ makes the triangle (III)
commute, and so the map we have constructed is a bijection.

\end{document}